\documentclass[numbook,envcountsame]{svjour3} 
\smartqed  
\usepackage{amssymb,amsmath,amsfonts}
\usepackage{enumitem}
\usepackage{hyperref}
\hypersetup{colorlinks=true,linkcolor=blue,urlcolor=blue,bookmarks=true}

\DeclareMathOperator{\haus}{dist_{\mathbb{R}}}
\DeclareMathOperator{\dist}{dist}
\DeclareMathOperator{\Dist}{Dist}
\DeclareRobustCommand{\epslabel}{\varepsilon}

\hypersetup{colorlinks=true,linkcolor=blue,urlcolor=blue,bookmarks=true}
\begin{document}

\title{Sigmoidal approximations of a nonautonomous neural network with infinite
  delay and Heaviside function
  \thanks{Both authors were partially supported by MICIIN/FEDER
    under project RTI2018-096523-B-100.} }


\author{Peter E. Kloeden \and V{\'\i}ctor M. Villarragut}


\dedication{Dedicated to the memory of Russell A. Johnson}

\institute{P.E. Kloeden \at
  Mathematics Department, University of T{\"u}bingen, T{\"u}bingen 72076, Germany \\
  \email{kloeden@na.uni-tuebingen.de} \and V.M. Villarragut \at Departamento de
  Matem{\'a}tica Aplicada a la Ingenier{\'\i}a Industrial, Universidad
  Polit{\'e}cnica de Madrid, Calle de Jos{\'e} Guti{\'e}rrez Abascal 2, 28006 Madrid, Spain \\
  \email{victor.munoz@upm.es} }

\date{}

\maketitle

\vspace*{-1.5cm}

\noindent\begin{tabular}{|l|}
  \hline
  Kloeden, P.E., Villarragut, V.M.\\Sigmoidal Approximations of a Nonautonomous
Neural Network with Infinite\\ Delay and Heaviside Function. J Dyn Diff Equat 34, 721–745 (2022).\\
https://doi.org/10.1007/s10884-020-09899-4\\
  \copyright ~Springer\\
  \hline
\end{tabular}

\vspace{1cm}

\begin{abstract}
  In this paper, we approximate a nonautonomous neural network with infinite
  delay and a Heaviside signal function by neural networks with sigmoidal
  signal functions. We show that the solutions of the sigmoidal models converge
  to those of the Heaviside inclusion as the sigmoidal parameter vanishes. In
  addition, we prove the existence of pullback attractors in both cases, and the convergence of
  the attractors of the sigmoidal models to those of the Heaviside inclusion.
\end{abstract}

\section{Introduction}

Neural networks arise in many contexts, ranging from biological modeling to
engineering applications. The study of the dynamics of neural networks has been
of paramount importance in the process of establishing a theoretical framework
for this fast-paced field. Extensive efforts have been devoted to the
development of such theory (see e.g. Wu~\cite{book:w} and the references
therein). Specifically, in~\cite{book:w}, the following system of differential
equations with atomic delay was proposed:
\[
  \begin{split}
    &x_i'(t)=-A_i(x_i(t))\,x_i(t) +\sum_{\substack{k=1\\ k\neq
        i}}^{n}(z_{ki}(t)-c_{ki})\,f_k(x_k(t-\tau_{ki})-\Gamma_k)+I_i(t)\,,\\
    &z_{ij}'(t)=-B_{ij}(z_{ij}(t))\,z_{ij}(t)
    +d_{ij}\,f_i(x_i(t-\tau_{ij})-\Gamma_i)\,[x_j(t)]^+\,,
  \end{split}
\]
where $i,j\in\{1,\ldots,n\}$ with $i\neq j$. The parameters in this model are
thoroughly defined in~\cite{book:w}, but we will give a few representative
definitions for the reader's convenience. Namely, $x_i(t)$ represents the
deviation of the $i$th neuron's potential from its equilibrium and $z_{ij}(t)$ is
the neurotransmitter average release rate per unit axon signal
frequency, usually referred to as synaptic coupling strength. We assume further that the neuron's potential decays exponentially to its
equilibrium in the absence of external processes, which is regulated by
$A_i$. Also, the $k$th neuron sends a signal $f_k$, which depends on the past
values of $x_k$ and is subject to a threshold $\Gamma_k$, along the axon to the
$i$th neuron. That signal affects the target neuron in an additive manner
proportionally to the coupling strength $z_{ki}(t)$, giving rise to the term $\sum_{k\neq
  i}z_{ki}(t)\,f_k(x_k(t-\tau_{ki})-\Gamma_k)$. Moreover, a hardwiring
of the inhibitory input from the $k$th neuron to the $i$th neuron is present,
that is, its coupling strength is constant, giving rise to the term $-\sum_{k\neq
  i}c_{ki}\,f_k(x_k(t-\tau_{ki})-\Gamma_k)$, with $c_{ki}\geq 0$. Finally, $I_i(t)$
represents the external stimuli arriving at the $i$th neuron. As for the second
equation, it is based on Hebb's law. Specifically, $d_{ij}$ are nonnegative
constants and $[z]^+=\max\{z,0\}$ for each $z\in\mathbb{R}$. It is noteworthy that the functions $f_i$ may be sigmoidal functions of the form
$\sigma_\varepsilon(x)=\frac{1}{1+e^{-x/\varepsilon}}$, $x\in\mathbb{R}$, or a
Heaviside function (cf. Levine~\cite{book:l}, and~McCulloch and
Pitts~\cite{paper:mp}).

More recently, Han and Kloeden~\cite{paper:hk1} addressed the study of a
nonautonomous neural network given by a lattice of differential inclusions and a
Heaviside signal function. In order to obtain a rigorous setting for these
inclusions, they replaced the Heaviside signal function with
\[
  s\in\mathbb{R} \mapsto\chi(s)=\left\{
    \begin{array}{ll}
      \{0\} & \text{if } s<0, \\ \relax 
      [0,1] & \text{if } s=0, \\ \relax
      \{1\} & \text{if } s>0.
    \end{array}
  \right.
\]
Specifically, they proved the existence of pullback attractors and forward
omega-limit sets. Han and Kloeden~\cite{paper:hk2} reached further conclusions,
proving the convergence of the solutions of a family of autonomous neural
networks given by a lattice of differential equations with sigmoidal signal
functions to solutions of the corres\-ponding differential inclusion with a
Heaviside signal function. They also proved the existence of pullback attractors
in these two cases and the convergence of the former to the latter in the
Hausdorff pseudometric. Further developments in this line can be found in Wang,
Kloeden, and Yang~\cite{paper:wky1} and~\cite{paper:wky2}, where the results of
the two previous papers were extended to the case of neural networks given by
autonomous differential equations and inclusions with atomic delay.

The introduction of distributed delays into models of neural networks dates back
to Herz, Salzer, K{\"u}hn, and van Hemmen~\cite{paper:hskv}. Differential
equations with infinite delay and the dynamics have been widely studied (see
Hino, Murakami, and Naito~\cite{book:hmn} for a thorough description of the
available phase spaces and some general results for this kind of equations).

An interesting generalization of the results on the dynamics of neural networks
presented above is to consider neural networks given by family of systems of
nonautonomous differential equations with sigmoidal signal function and infinite
delay, and the corresponding system of inclusions with a Heaviside signal
function and infinite delay. Such systems are given, for each
$i,j\in\{1,\ldots,n\}$ with $i\neq j$, as follows:
\[
  \begin{split}
    & x_i'(t)=-A_i(x_i(t))\,x_i(t)+\sum_{k\neq
      i}(z_{ki}(t)-c_{ki})\int_{-\infty}^{0}
    \!\!\sigma_\varepsilon(x_k(t+\tau)-\Gamma_k)\,d\mu_{ki}(\tau)+I_i(t),\\
    & z_{ij}'(t)=-B_{ij}(z_{ij}(t))\,z_{ij}(t)
    +d_{ij}\,\int_{-\infty}^{0}\sigma_\varepsilon
    (x_i(t+\tau)-\Gamma_i)\,d\mu_{ij}(\tau)\,[x_j(t)]^+\,,
  \end{split}
\]
together with the corresponding system of differential inclusions obtained by
replacing $\sigma_\varepsilon$ with $\chi$ and where the integrals must be
understood in the Aumann sense, as detailed
in~Section~\ref{sec:preliminaries}. The parameters in this model are as above,
except for the measures $\mu_{ij}$, which allow for much more general delay
terms.

The goal of this work is to provide approximations of the solutions of the
aforementioned inclusion by those of the family of equations, to prove the
existence of pullback attractors in both cases, and to guarantee the convergence
of the pullback attractors of the family of sigmoidal equations to that of the
Heaviside inclusion in the Hausdorff pseudometric.

In Section~\ref{sec:preliminaries}, we introduce some necessary notation, as
well as some useful definitions and properties of the Hausdorff pseudometric and
the Aumann integral. Besides, we recall some previous results on the asymptotic
behavior of some set-valued processes, which will be used in
Section~\ref{sec:asymptotic}. Section~\ref{sec:existence} is devoted to the
existence and uniqueness theory for the proposed family of differential
equations with sigmoidal signal functions. Section~\ref{sec:existence} also
addresses the convergence of the solutions of the family of differential
equations with sigmoidal signal functions to those of the differential inclusion
with Heaviside signal function. In Section~\ref{sec:asymptotic}, we establish
the existence of a unique global attractor for the equations of the family of
differential equations and for the differential inclusion. Finally, in
Section~\ref{sec:comparison}, we prove the convergence of the pullback
attractors of the family of sigmoidal equations to that of the Heaviside
inclusion in the Hausdorff pseudometric.

\section{Preliminaries}\label{sec:preliminaries}

\subsection{Notation}

We shall identify $\mathbb{R}^{n(n-1)}$ with the set
\[
  \{(z_{ij})_{i,j\in\{1,\ldots,n\},\,i\neq j}:
  z_{ij}\in\mathbb{R},\,i,j=1,\ldots,n,\,i\neq j\}
\]
and denote its elements by $(z_{ij})_{i\neq j}$. Likewise, we shall identify
$\mathbb{R}^{n^2}$ with the set $\mathbb{R}^n\times \mathbb{R}^{n(n-1)}$ and
denote its elements by $(x,z)=((x_i)_i,(z_{ij})_{i\neq j})$. Let
$\pi_i:\mathbb{R}^{n^2}\to\mathbb{R}$ and
$\pi_{ij}:\mathbb{R}^{n^2}\to\mathbb{R}$ denote the $n^2$ canonical projections
from $\mathbb{R}^{n^2}$ onto $\mathbb{R}$. We will abuse notation by letting
$\left\|{\cdot}\right\|$ denote the supremum norm on both $\mathbb{R}^n$ and
$\mathbb{R}^{n(n-1)}$, as well as the following norm on $\mathbb{R}^{n^2}$:
\[
  \left\|u\right\|=\max\left\{\left\|x\right\|,\left\|z\right\|\right\},\quad
  u=(x,z)\in\mathbb{R}^{n^2}\,.
\]

Throughout this work, $C_\gamma$ will denote the space
$C_\gamma((-\infty,0],\mathbb{R}^{n^2})$, for some $\gamma>0$, endowed with the
norm $\left\|u\right\|_\gamma=\sup_{t\leq 0}e^{\gamma\,t}\left\|u(t)\right\|$,
$u\in C_\gamma$. Given $r>0$, we will also write
$\overline{B}_\gamma(0,r)=\{u\in C_\gamma:\left\|u\right\|_\gamma\leq
r\}$. Moreover, given $x:(-\infty,a]\to\mathbb{R}^m$ for some $a\in\mathbb{R}$
and some $m\in\mathbb{N}$, we will denote its \emph{history} up to
$t\in(-\infty,a]$ by $x_t:(-\infty,0]\to\mathbb{R}^m$, $s \mapsto x(t+s)$. See
e.g.~\cite{book:hmn} for the main properties of $C_\gamma$.

The following result will be useful throughout this work.

\begin{lemma}\label{gamma_implies_compact_open}
  If $\{u^m\}_m\subset C_\gamma$ converges to $u\in C_\gamma$, then $\{u^m\}_m$
  converges to $u$ uniformly on compact sets.
\end{lemma}
\begin{proof}
  Fix $\rho<0$. Then
  $\sup_{\rho\leq s\leq 0}\left\|u^m(s)-u(s)\right\| \leq
  e^{-\gamma\,\rho}\left\|u^m-u\right\|_\gamma$. The result follows immediately.
  \qed
\end{proof}

We will also consider the Banach space of Borel signed measures on $(-\infty,0]$
endowed with the total variation norm, which we will denote by $|{\cdot}|$.

Finally, if $f:A\to X$ is a mapping into a normed space $(X,\|{\cdot}\|_X)$ and
$S\subset A$, then we will write
$\left\|f\right\|_S=\sup_{a\in S}\left\|f(a)\right\|_X$.

\subsection{Set-valued topology and integrals}

In this work, we will make extensive use of the Hausdorff pseudometric on
$\mathbb{R}$.

\begin{definition}
  Let $A,\,B$ be nonempty subsets of $\mathbb{R}$. Then
  \[
    \haus(A,B)=\sup_{a\in A}\inf_{b\in B} |a-b| =\inf\{\varepsilon>0:A\subset
    B+\varepsilon\,[-1,1]\}
  \]
  is a pseudometric called the \emph{Hausdorff pseudometric}.
\end{definition}

The following lemma provides a useful property of the Hausdorff
pseudometric. Its proof is included in that of Lemma~2.2 of De
Blasi~\cite{paper:d}. See also Proposition~2.4.1 of Diamond and
Kloeden~\cite{book:dk} for further details.

\begin{lemma}\label{properties_hausdorff}
  If $A$, $A_1$, $B$, $B_1$ are nonempty and bounded subsets of $\mathbb{R}$,
  then
  \begin{enumerate}[label={\upshape(\roman*)}]
  \item $\haus(t\,A,t\,B)=t\,\haus(A,B)$ for all $t\geq0$;
  \item $\haus(A+B,A_1+B_1)\leq\haus(A,A_1)+\haus(B,B_1)$.
  \end{enumerate}
\end{lemma}

Finally, the following definition, in the line of Definition~2.7 of Pucci and
Viti\-llaro~\cite{paper:pv}, will allow us to study some differential inclusions
with infinite delay.

\begin{definition}
  Let $F:(-\infty,0]\to\mathcal{P}(\mathbb{R})$ be a measurable set-valued
  function. The \emph{Aumann integral} of $F$ with respect to a Borel measure
  $\mu$ is
    \[
      \int_{-\infty}^{0}F\,d\mu=\bigg\{\int_{-\infty}^{0}s\,d\mu:
      \begin{array}{l}
        s \text{ is
          a selector of }F \text{ and }s \text{ is integrable} \\
          \text{with respect to }\mu
      \end{array}
      \bigg\}.
  \]
\end{definition}
Some properties of the Aumann integral in the setting of the present work will
be proved in Section~\ref{sec:existence}.

\subsection{Set-valued dynamical systems}

Let $X$ be a Banach space with norm $\left\|{\cdot}\right\|_X$ and let
$\mathcal{P}_c(X)$ denote the family of nonempty closed subsets of $X$. Let
$\mathbb{R}^2_\geq=\{(t,t_0)\in\mathbb{R}^2:t\geq t_0\}$.

\begin{definition}
  A mapping $\Phi:\mathbb{R}^2_\geq\times X\to\mathcal{P}_c(X)$ is said to be a
  \emph{set-valued process} if the following conditions hold:
  \begin{enumerate}[label={\upshape(\roman*)}]
  \item $\Phi(t_0,t_0,x)=\{x\}$ for all $x\in X$ and all $t_0\in\mathbb{R}$;
  \item $\Phi(t_2,t_0,x)=\Phi(t_2,t_1,\Phi(t_1,t_0,x))$ for all $x\in X$ and all
    $t_0,t_1,t_2\in\mathbb{R}$ with \mbox{$t_2\geq t_1\geq t_0$}.
  \end{enumerate}
\end{definition}

Let $\mathcal{M}(X)$ denote the family of nonautonomous subsets of $X$ whose
images are nonempty and closed.

\begin{definition}
  Let $\Phi$ be a set-valued process on $X$ and let
  $\mathcal{D}\subset\mathcal{M}(X)$. Let $\{\Lambda(t)\}_{t\in\mathbb{R}}$ be a
  nonautonomous subset of $X$.
  \begin{enumerate}[label={\upshape(\roman*)}]
  \item $\{\Lambda(t)\}_{t\in\mathbb{R}}$ is \emph{positively invariant}
    (resp. \emph{invariant}) for $\Phi$ if
    $\Phi(t,t_0,\Lambda(t_0))\subset\Lambda(t)$
    (resp. $\Phi(t,t_0,\Lambda(t_0))=\Lambda(t)$) for all
    $(t,t_0)\in\mathbb{R}^2_\geq$.
  \item $\{\Lambda(t)\}_{t\in\mathbb{R}}$ is \emph{$\mathcal{D}$-pullback
      absorbing} for $\Phi$ if, for each
    $\{D(t)\}_{t\in\mathbb{R}}\in\mathcal{D}$ and each $t\in\mathbb{R}$, there
    is a $T>0$ such that
    \[
      \Phi(t,t-\tau,D(t-\tau))\subset\Lambda(t)\quad \text{for all }\tau\geq
      T\,.
    \]
  \item $\{\Lambda(t)\}_{t\in\mathbb{R}}$ is \emph{$\mathcal{D}$-pullback
      attracting} for $\Phi$ if, for each
    $\{D(t)\}_{t\in\mathbb{R}}\in\mathcal{D}$,
    \[
      \lim_{\tau\to \infty}\dist_X(\Phi(t,t-\tau,D(t-\tau)),\,\Lambda(t))=0\,,
    \]
    where $\dist_X(A,B)=\sup_{a\in A}\inf_{b\in B}\left\|a-b\right\|_X$.
  \end{enumerate}
\end{definition}

\begin{definition}
  A subset $\mathcal{D}\subset\mathcal{M}(X)$ is said to be a \emph{universe} if
  it is closed for inclusions, i.e. if $\{D(t)\}_{t\in\mathbb{R}}\in\mathcal{D}$
  and $\{\widetilde{D}(t)\}_{t\in\mathbb{R}}\in\mathcal{M}(X)$ with
  $\widetilde{D}(t)\subset D(t)$ for all $t\in\mathbb{R}$, then
  $\{\widetilde{D}(t)\}_{t\in\mathbb{R}}\in\mathcal{D}$.
\end{definition}

\begin{definition}
  Let $\Phi$ be a set-valued process on $X$ and let $\mathcal{D}$ be a universe.
  A nonautonomous subset $\mathcal{A}=\{A(t)\}_{t\in\mathbb{R}}$ of $X$ is a
  \emph{global $\mathcal{D}$-pullback attractor} for $\Phi$ if
  \begin{enumerate}[label={\upshape(\roman*)}]
  \item $A(t)$ is a compact subset of $X$ for all $t\in\mathbb{R}$,
  \item $\mathcal{A}$ is invariant, and
  \item $\mathcal{A}$ is $\mathcal{D}$-pullback attracting.
  \end{enumerate}
\end{definition}

The existence of a pullback attractor relies on the continuity properties and
the asymptotic compactness of a set-valued process. The following definitions
will be used in order to prove such existence.

\begin{definition}
  A set-valued process $\Phi$ on $X$ is \emph{upper semi continuous at
    $x_0\in X$} if, for each $(t,t_0)\in\mathbb{R}^2_\geq$ and neighborhood $V$
  of the subset $\Phi(t,t_0,x_0)$ of $X$, there is a $\delta>0$ such that
  $\Phi(t,t_0,x)\in V$ for all $x\in X$ with $\left\|x-x_0\right\|_X<\delta$. A
  set-valued process $\Phi$ is \emph{upper semi continuous} if it is upper semi
  continuous at every $x_0\in X$.
\end{definition}

\begin{definition}
  Let $\Phi$ be a set-valued process on $X$ and let $\mathcal{D}$ be a
  universe. $\Phi$ is said to be \emph{$\mathcal{D}$-pullback asymptotically
    upper semi compact} in $X$ if the following property holds: if
  $t\in\mathbb{R}$, $\tau_m\to\infty$ as $m\to\infty$,
  $\{D(t)\}_{t\in\mathbb{R}}\in\mathcal{D}$, and \mbox{$\{x_m\}_m\subset X$} with
  $x_m\in D(t-\tau_m)$ for all $m\in\mathbb{N}$, then every sequence
  $\{y_m\}_m\subset X$ with \mbox{$y_m\in\Phi(t,t-\tau_m,x_m)$} for all
  $m\in\mathbb{N}$ has a convergent subsequence.
\end{definition}

The following result, presented in Caraballo and Kloeden~\cite{paper:ck},
guarantees the existence of pullback attractors. It will be used in the proofs
of the main results in Section~\ref{sec:asymptotic}.

\begin{theorem}\label{abstract_attractor}
  Let $\Phi$ be a set-valued process on $X$ and let $\mathcal{D}$ be a
  universe. Assume that
  \begin{enumerate}[label={\upshape(\roman*)}]
  \item $\Phi$ is upper semi continuous,
  \item $\Phi$ has a $\mathcal{D}$-pullback absorbing set
    $\{\Lambda(t)\}_{t\in\mathbb{R}}\in\mathcal{D}$,
  \item $\Phi$ is $\mathcal{D}$-pullback asymptotically upper semi compact in
    $X$.
  \end{enumerate}
  Then $\Phi$ has a $\mathcal{D}$-pullback attractor $\{A(t)\}_{t\in\mathbb{R}}$
  given by
  \[
    A(t)=\bigcap_{t_0\geq 0}\overline{\bigcup_{\tau\geq
        t_0}\Phi(t,t-\tau,\Lambda(t-\tau))}\,, \qquad t\in\mathbb{R}\,.
  \]
\end{theorem}

\section{Existence and uniqueness of solutions}\label{sec:existence}

\subsection{The sigmoidal systems}

Let us consider the family of neural networks \addtocounter{equation}{1}
\begin{equation}\label{eq:sigmoidal}\tag*{{\upshape(\theequation)$_{\epslabel}$}}
  (x,z)'(t)=D(x_t,z_t)+\Sigma^\varepsilon(x_t,z_t)+I(t),\quad t\geq t_0\,,
\end{equation}
where $t_0\in\mathbb{R}$, $\varepsilon\in(0,1)$, and $D$, $\Sigma^\varepsilon$
and $I$ are defined as follows. $D$ is the \emph{decay} term and it is defined
by
\[
  \begin{array}[t]{rccl}
    D: & C_\gamma & \longrightarrow & \mathbb{R}^{n^2} \\
       & (x,z) & \mapsto         &
                                   \big((-A_i(x_i(0))\,x_i(0))_i,\,(-B_{ij}(z_{ij}(0))\,z_{ij}(0))_{i\neq j}\big)\,,
  \end{array}
\]
where $A_i,\,B_{ij}:\mathbb{R}\to\mathbb{R}$ for all $i,j\in\{1,\ldots,n\}$ with
$i\neq j$. Notice that, in this case, $(x,z):(-\infty,0]\to\mathbb{R}^{n^2}$,
and its components are given by $x_i=\pi_i\circ (x,z)$
 and $z_{ij}=\pi_{ij}\circ (x,z)$ for all $i,j\in\{1,\ldots,n\}$ with
$i\neq j$.

For each $\varepsilon\in(0,1)$, we consider the \emph{sigmoidal
  function}
\[
  \begin{array}[t]{rccl}
    \sigma_\varepsilon: & \mathbb{R} & \longrightarrow & \mathbb{R} \\
                        & x & \mapsto & \frac{1}{1+e^{-x/\varepsilon}}\,.
  \end{array}
\]
This is a globally Lipschitz function with Lipschitz constant
$\frac{1}{\varepsilon}$. $\Sigma^\varepsilon$ is the \emph{excitation and
  inhibition} term and it is defined by
\[
  \begin{array}[t]{rccl}
    \Sigma^\varepsilon: & C_\gamma & \longrightarrow & \mathbb{R}^{n^2} \\
                        & (x,z) & \mapsto         & \Big( \left(\sum_{k\neq
                                                    i}(z_{ki}(0)-c_{ki})\,\int_{-\infty}^{0}\sigma_\varepsilon
                                                    (x_k(\tau)-\Gamma_k)\,d\mu_{ki}(\tau)\right)_i,\\
                        &&&\quad \left(d_{ij}\,\int_{-\infty}^{0}\sigma_\varepsilon
                            (x_i(\tau)-\Gamma_i)\,d\mu_{ij}(\tau)\,[x_j(0)]^+\right)_{i\neq j}\Big)\,,
  \end{array}
\]
where $\mu_{ij}$ is a positive regular Borel measure on $(-\infty,0]$,
$c_{ij},\,d_{ij}\geq 0$ for all $i,j\in\{1,\ldots,n\}$ with $i\neq j$, and
$[{\cdot}]^+:\mathbb{R}\to\mathbb{R}$, $t\mapsto\max\{t,0\}$. We will write
$c=(c_{ij})_{i\neq j}$ and $d=(d_{ij})_{i\neq j}$ for simplicity. Finally, $I$
is the \emph{stimuli} term and it is defined by
\[
  \begin{array}[t]{rccl}
    I: & \mathbb{R} & \longrightarrow & \mathbb{R}^{n^2} \\
       & t & \mapsto         & \left((I_i(t))_i,0\right)\,,
  \end{array}
\]
where $I_i:\mathbb{R}\to\mathbb{R}$ is a given function for all
$i\in\{1,\ldots,n\}$.

For each $\varepsilon\in(0,1)$, let us define
\[
  \begin{array}[t]{rccl}
    f_\varepsilon: & \mathbb{R}\times C_\gamma & \longrightarrow & \mathbb{R}^{n^2} \\
                   & (t,x,z) & \mapsto         & D(x,z)+\Sigma^\varepsilon(x,z)+I(t)\,.
  \end{array}
\]
Then equation~\ref{eq:sigmoidal} can be written as
\[
  (x,z)'(t)=f_\varepsilon(t,x_t,z_t),\quad t\geq t_0\,.
\]

We will assume the following hypotheses:

\begin{enumerate}[label={\upshape (D)}]
\item\label{hyp:D} $A_i$ and $B_{ij}$ are locally Lipschitz continuous on
  $\mathbb{R}$ and there exists $\alpha>0$ such that $A_i(s)\geq\alpha$ and
  $B_{ij}(s)\geq\alpha$ for all $s\in\mathbb{R}$ and all $i,j\in\{1,\ldots,n\}$
  with $i\neq j$.
\end{enumerate}
\begin{enumerate}[label={\upshape (M)}]
\item\label{hyp:M} the measure $\mu_{ij}$ satisfies
  $\int_{-\infty}^{0}e^{-\gamma\,t}d\mu_{ij}(t)<\infty$ for all
  $i,j\in\{1,\ldots,n\}$ with $i\neq j$.
\end{enumerate}
\begin{enumerate}[label={\upshape (I)}]
\item\label{hyp:I} $I_i$ is continuous on $\mathbb{R}$ for all
  $i\in\{1,\ldots,n\}$.
\end{enumerate}

Notice that, in particular,~\ref{hyp:D} implies that $A_i$ and $B_{ij}$ are
Lipschitz continuous on all the intervals $[-r,r]$, $r>0$.

We present some results regarding the main properties of the function
$f_\varepsilon$.

\begin{proposition}\label{f_continuous}
  The function $f_\varepsilon$ is continuous for all $\varepsilon>0$.
\end{proposition}
\begin{proof}
  The continuity of $f_\varepsilon$ in its first variable follows
  from~\ref{hyp:I}. It remains to prove that $f_\varepsilon$ is continuous in
  its second variable. Indeed, fix $\{u^k\}_k\subset C_\gamma$ such that
  \mbox{$u^k=(x^k,z^k)\to u=(x,z)\in C_\gamma$} as $k\to\infty$. Then
  $u^k(0)\to u(0)$, and, consequently, $D(u^k)\to D(u)$ as $k\to\infty$ thanks
  to~\ref{hyp:D}.

  On the other hand, for all $i,j\in\{1,\ldots,n\}$ with $i\neq j$,
  \begin{multline*}
    \bigg|\int_{-\infty}^{0}\sigma_\varepsilon(x^k_i(\tau)-\Gamma_i)d\mu_{ij}(\tau)
    -\int_{-\infty}^{0}\sigma_\varepsilon(x_i(\tau)-\Gamma_i)d\mu_{ij}(\tau)\bigg|\leq\\
    \leq \frac{1}{\varepsilon}\,
    \int_{-\infty}^{0}\left|x^k_i(\tau)-x_i(\tau)\right|d\mu_{ij}(\tau) \leq
    \frac{1}{\varepsilon}\,
    \int_{-\infty}^{0}e^{-\gamma\,\tau}d\mu_{ij}(\tau)\,\left\|x^k-x\right\|_\gamma\,,
  \end{multline*}
  which vanishes as $k\to\infty$ thanks to~\ref{hyp:M}. Finally, the continuity
  of $[{\cdot}]^+$ yields $\Sigma^\varepsilon(u^k)\to\Sigma^\varepsilon(u)$, as
  wanted.  \qed
\end{proof}

\begin{proposition}\label{f_bounded}
  For all $a,b\in\mathbb{R}$ with $a\leq b$ and all $r>0$, there exists $L>0$
  such that $f_\varepsilon$ is bounded by $L$ on
  $[a,b]\times \overline{B}_\gamma(0,r)$ for all $\varepsilon>0$.
\end{proposition}
\begin{proof}
  Hypothesis~\ref{hyp:I} and the compactness of $[a,b]$ yield the boundedness of
  $I$ on $[a,b]$. Fix $u=(x,z)\in \overline{B}_\gamma(0,r)$. Then
  $\left\|x(0)\right\|\leq \left\|u(0)\right\|\leq \left\|u\right\|_\gamma\leq
  r$ and, analogously, $\left\|z(0)\right\|\leq r$. Hence, for all
  $i,j\in\{1,\ldots,n\}$ with $i\neq j$,
  \[
    \left|-A_i(x_i(0))\,x_i(0)\right|\leq \sup_{s\in[-r,r]}A_i(s)\,r\,,\quad
    \left|-B_{ij}(z_{ij}(0))\,z_{ij}(0)\right|\leq
    \sup_{s\in[-r,r]}B_{ij}(s)\,r\,.
  \]
  This, together with~\ref{hyp:D}, implies that $D$ is bounded. Notice that
  $0\leq\sigma_\varepsilon\leq 1$ for all $\varepsilon>0$,
  $[x_j(0)]^+\leq |x_j(0)|\leq r$ and $|z_{ij}(0)|\leq r$, whence
  $\Sigma^\varepsilon$ is bounded thanks to~\ref{hyp:M}, as desired.  \qed
\end{proof}

In the present setting, we are able to prove the local existence of solutions of
the sigmoidal system of differential equations~\ref{eq:sigmoidal}.

\begin{theorem}
  Assume hypotheses~\ref{hyp:D},~\ref{hyp:M}, and~\ref{hyp:I}. Then, for each
  $\varepsilon>0$, $t_0\in\mathbb{R}$ and $u_0\in C_\gamma$, the system of
  equations~\ref{eq:sigmoidal} admits a local solution defined on
  $[t_0,t_1)$, $t_1>t_0$, with initial datum $u_0$ at $t=t_0$.
\end{theorem}
\begin{proof}
  Following Theorem 3.7 on p.15 and Theorem 1.1 on p.36 of~\cite{book:hmn}, we
  define $K:\mathbb{R}\to\mathbb{R}$, $t \mapsto\max\{1,e^{-\gamma\,t}\}$. Then,
  for all $m>0$,
  \[
    \lim_{t\to 0^+}K(t)\,\int_{0}^{t}m\,ds=\lim_{t\to
      0^+}\max\{m\,t,\,m\,t\,e^{-\gamma\,t}\}=0\,.
  \]
  As a result, there is a $t_1>0$ such that the conditions of the
  aforementioned theorems hold. This, together with
  Propositions~\ref{f_continuous} and~\ref{f_bounded}, which guarantee the
  Caratheodory character of $f_\varepsilon$, yields the result.  \qed
\end{proof}

\begin{lemma}\label{history_continuous}
  If $a\in\mathbb{R}$, $u:(-\infty,a)\to\mathbb{R}^{n^2}$ is continuous and
  $u_{t_0}\in C_\gamma$ for some $t_0<a$, then the mapping
  $(-\infty,a)\to C_\gamma$, $t \mapsto u_t$ is continuous.
\end{lemma}
\begin{proof}
  It follows from Theorem 3.7 on p.15 of~\cite{book:hmn}.  \qed
\end{proof}

Let us prove one more property as to the Lipschitz character of $f_\varepsilon$.

\begin{proposition}\label{f_Lipschitz}
  For all $\varepsilon>0$, $f_\varepsilon$ is Lipschitz continuous in its spatial
  variables on sets of the form $\mathbb{R}\times \overline{B}_\gamma(0,r)$,
  $r>0$.
\end{proposition}
\begin{proof}
  Fix $t\in\mathbb{R}$ and
  $u^1=(x^1,z^1),u^2=(x^2,z^2)\in \overline{B}_\gamma(0,r)$. For all
  $i,j\in\{1,\ldots,n\}$ with $i\neq j$, clearly,
  \[
    \left|\int_{-\infty}^{0}(\sigma_\varepsilon(x^1_i(\tau)-\Gamma_i)
      -\sigma_\varepsilon(x^2_i(\tau)-\Gamma_i))d\mu_{ij}(\tau)\right|
    \leq\frac{1}{\varepsilon}\int_{-\infty}^{0}e^{-\gamma\,\tau}
    d\mu_{ij}(\tau)\,\left\|u^1-u^2\right\|_\gamma\,.
  \]
  Besides, $|z_{ij}(0)|\leq r$. As a result, \ref{hyp:M} implies that there is
  an $L_1>0$ independent of $t$, $u^1$ and $u^2$ such that
  \begin{equation}
    \label{eq:Lipschitz_Sigma}
    \left\|\Sigma^\varepsilon(u^1)-\Sigma^\varepsilon(u^2)\right\|\leq
    L_1\,\left\|u^1-u^2\right\|_\gamma\,.  
  \end{equation}

  On the other hand, fix $i,j\in\{1,\ldots,n\}$ with $i\neq
  j$. From~\ref{hyp:D}, it follows that $A_i$, $B_{ij}$ are bounded by $K>0$ on
  $[-r,r]$ and they have a common Lipschitz constant $L>0$ on $[-r,r]$. Then
  \begin{equation}
    \label{eq:Lipschitz_D}
    \begin{split}
      |A_i(x^1_i&(0))\,x^1_i(0)-A_i(x^2_i(0))\,x^2_i(0)|\leq\\ &\leq
      K\,|x^1_i(0)-x^2_i(0)|+r\,L\,|x^1_i(0)-x^2_i(0)|
      \leq(K+r\,L)\left\|u^1-u^2\right\|_\gamma, \\
      |B_{ij}(z^1_{ij}&(0))\,z^1_{ij}(0)-B_{ij}(z^2_{ij}(0))\,z^2_{ij}(0)|
      \leq(K+r\,L)\left\|u^1-u^2\right\|_\gamma\,.
    \end{split}
  \end{equation}
  Finally, from~\eqref{eq:Lipschitz_Sigma}, \eqref{eq:Lipschitz_D} and the fact
  that $I(t)-I(t)=0$, it follows that $f_\varepsilon$ satisfies the desired
  property.  \qed
\end{proof}

The preceding result allows us to guarantee the uniqueness of the solutions of
the sigmoidal system~\ref{eq:sigmoidal}.

\begin{theorem}\label{sigmoidal_uniqueness}
  Under hypotheses~\ref{hyp:D},~\ref{hyp:M}, and~\ref{hyp:I}, for each
  $\varepsilon>0$, $t_0\in\mathbb{R}$ and $u_0\in C_\gamma$, the system of
  equations~\ref{eq:sigmoidal} admits at most one solution with initial datum
  $u_0$ at $t=t_0$.
\end{theorem}
\begin{proof}
  It follows from Proposition~\ref{f_Lipschitz} and Theorem 1.2 on p.38
  of~\cite{book:hmn}.  \qed
\end{proof}

The following result provides a bound for the solutions of the sigmoidal
system~\ref{eq:sigmoidal}.

\begin{proposition}\label{general_bound}
  Fix $[a,b]\subset\mathbb{R}$, $r>0$ and $T>0$. There exist $k_1$, $k_2>0$ such
  that, for all $\varepsilon>0$ and all $u:(-\infty,t_0+T)\to\mathbb{R}^{n^2}$
  solution of~\ref{eq:sigmoidal} with $t_0\in[a,b]$ and initial datum
  $u_{t_0}\in \overline{B}_\gamma(0,r)$, $u$ satisfies
  \[
    \left\|u_t\right\|_\gamma\leq k_1\,e^{k_2\,(t-t_0)}\quad \text{for all
    }t\in[t_0,t_0+T)\,.
  \]
\end{proposition}
\begin{proof}
  Fix $t\in[t_0,t_0+T)$ and let $u=(x,z)$. Then
  \begin{equation}
    \label{eq:solution_bound_history}
    \left\|u_t\right\|_\gamma=\max\left\{\sup_{s\leq
        t_0-t}e^{\gamma\,s}\left\|u(t+s)\right\|,
      \,\sup_{t_0-t\leq s\leq 0}e^{\gamma\,s}\left\|u(t+s)\right\|\right\}\,.
  \end{equation}
  On the one hand,
  $\sup_{s\leq
    t_0-t}e^{\gamma\,s}\left\|u(t+s)\right\|=e^{-\gamma(t-t_0)}\left\|u_{t_0}\right\|_\gamma\leq
  r$. On the other hand, thanks to~\ref{hyp:D}, for all $i,j\in\{1,\ldots,n\}$
  with $i\neq j$, we have
  \[
    -A_i(x_i(t))\,x_i^2(t)\leq-\alpha\,x_i^2(t)\leq-\alpha\,\|u(t)\|_2^2 \quad \text{and}\quad
    -B_{ij}(z_{ij}(t))\,z_{ij}^2(t)\leq-\alpha\,\|u(t)\|_2^2\,.
  \]
  Besides,
  \begin{multline*}
    \pi_i\circ\Sigma^\varepsilon(u_t)\,x_i(t)\leq
    \left|\pi_i\circ\Sigma^\varepsilon(u_t)\right|\,\left|x_i(t)\right|
    \leq \sum_{k\neq i}\left(|z_{ki}(t)|+c_{ki}\right)\,|\mu_{ki}|\,|x_i(t)|\\
    \leq\sum_{k\neq i}\left(\frac{z_{ki}^2(t)\,|\mu_{ij}|^2}{2}+
      \frac{x_i^2(t)}{2}\right)+\sum_{k\neq
      i}\left(\frac{c_{ki}^2\,|\mu_{ij}|^2}{2}+ \frac{x_i^2(t)}{2}\right)\leq
    k_1\,\left\|u(t)\right\|_2^2+k_2
  \end{multline*}
  for some $k_1$, $k_2>0$, where we have applied Young's inequality. Likewise,
  \[
    \pi_{ij}\circ\Sigma^\varepsilon(u_t)\,z_{ij}(t)\leq
    d_{ij}\,|\mu_{ij}|\,|x_j(t)|\,|z_{ij}(t)| \leq
    \frac{d_{ij}^2\,|\mu_{ij}|^2\,z_{ij}^2(t)}{2}+ \frac{x_j^2(t)}{2} \leq
    k_3\,\left\|u(t)\right\|_2^2
  \]
  for some $k_3>0$. Finally, there is a $k_4>0$ bound of $I_i$ on $[a,b+T]$
  thanks to~\ref{hyp:I}, whence
  \[
    I_i(t)\,x_i(t)\leq |I_i(t)|\,|x_i(t)|\leq k_4\,|x_i(t)|\leq
    \frac{k_4^2}{2}+\frac{x_i^2(t)}{2} \leq k_5\,\left\|u(t)\right\|_2^2+k_6
  \]
  for some $k_5$, $k_6>0$, where again we have used Young's
  inequality. Altogether, for all $i,j\in\{1,\ldots,n\}$ with $i\neq j$, we have
  \[
    \begin{split}
      &\tfrac{1}{2}(x_i^2(t))'=x'_i(t)\,x_i(t)\leq
      (-\alpha+k_1+k_5)\,\left\|u(t)\right\|_2^2+k_2+k_6 \quad \text{and}\\
      &\tfrac{1}{2}(z_{ij}^2(t))'=z'_{ij}(t)\,z_{ij}(t)\leq
      (-\alpha+k_3)\,\left\|u(t)\right\|_2^2\,.
    \end{split}
  \]
  By adding all inequalities, we conclude that there is a $k>0$ such that
  \[
    \frac{d}{dt}\left\|u(t)\right\|_2^2\leq k\,\left\|u(t)\right\|_2^2+k\,,\quad
    t\in[t_0,t_0+T)\,.
  \]
  Now, multiplying by $e^{-k\,t}$ on both sides and integrating, we have
  \[
    \left\|u(t)\right\|^2\leq\left\|u(t)\right\|_2^2\leq
    e^{k(t-t_0)}\left(1+\left\|u(t_0)\right\|_2^2\right)-1 \leq
    e^{k(t-t_0)}\left(1+n^2\,r^2\right).
  \]
  As a result,
  $\sup_{t_0-t\leq s\leq 0}e^{\gamma\,s}\left\|u(t+s)\right\|\leq
  \sqrt{1+n^2\,r^2}\,e^{\frac{k}{2}(t-t_0)}$. This, together with
  \eqref{eq:solution_bound_history}, yields the desired result.  \qed
\end{proof}

Thanks to the previous result, we are able to prove that the solutions of the
sigmoidal system~\ref{eq:sigmoidal} are, in fact, defined globally.

\begin{theorem}\label{global_existence}
  Assume hypotheses~\ref{hyp:D},~\ref{hyp:M}, and~\ref{hyp:I}. For each
  $\varepsilon>0$, the solutions of the system of equations~\ref{eq:sigmoidal}
  are defined globally.
\end{theorem}
\begin{proof}
  Fix $t_0\in\mathbb{R}$ such that the system of equations~\ref{eq:sigmoidal} is
  posed for all $t\geq t_0$. Suppose that $u:(-\infty,t_0+T)\to\mathbb{R}^{n^2}$
  is a noncontinuable solution of~\ref{eq:sigmoidal} for some $T>0$. Then,
  thanks to Proposition~\ref{general_bound}, there exist $k_1$, $k_2>0$ such
  that $\left\|u_t\right\|_\gamma\leq k_1\,e^{k_2\,(t-t_0)}$ for all
  $t\in[t_0,t_0+T)$. On the other hand, Theorem~2.7 on p.49 of~\cite{book:hmn}
  guarantees the existence of $\{t_k\}_k\subset[t_0,t_0+T)$ such that
  $\lim_{k\to \infty}t_k=t_0+T$ and
  $\left\|u_{t_k}\right\|_\gamma>k_1\,e^{k_2\,T}$, a contradiction. The result
  follows immediately.  \qed
\end{proof}

\subsection{The differential inclusions}

Let us turn to the study of some differential inclusions related to the
sigmoidal systems~\ref{eq:sigmoidal}. For all $\varepsilon>0$, let
$\chi_\varepsilon:\mathbb{R}\to\mathcal{P}(\mathbb{R})$ be defined for each
$s\in\mathbb{R}$ by
\[
  \chi_\varepsilon(s)=\left\{
    \begin{array}{ll}
      [0,\varepsilon] & \text{if } s<-b(\varepsilon), \\ \relax 
      [0,1] & \text{if } -b(\varepsilon)\leq s\leq b(\varepsilon), \\ \relax
      [1-\varepsilon,1] & \text{if } s>b(\varepsilon),
    \end{array}
  \right.
\]
where $b:[0,1]\to\mathbb{R}$ is given by
\[
  b(\varepsilon)=\left\{
    \begin{array}{ll}
      \varepsilon\,\log\left(\frac{1-\varepsilon}{\varepsilon}\right) & \text{if } \varepsilon\in(0,1), \\
      0 & \text{if } \varepsilon=0.
    \end{array}
  \right.
\]
It is noteworthy that $b$ is strictly increasing on, at least, the interval
$\left[0,\frac{1}{5}\right]$. We will denote $\chi=\chi_0$. Note that $\chi$ and
$\chi_\varepsilon$ are measurable with respect to the Borel $\sigma$-algebra
(see e.g. Definition~2.1 of~\cite{paper:pv} for further details).

We present some lemmas on the properties of the Aumann integral. In
particu\-lar, the following lemma circumvents the limitations of Theorem~3.1
of~\cite{paper:pv} and Theorem~8.6.3 of Aubin and Frankowska~\cite{book:af},
extending these results to --possibly-- atomic measures in this very specific
setting.

\begin{lemma}\label{convex_integral}
  Let $x:(-\infty,0]\to\mathbb{R}$ be a measurable function and let $\mu$ be a
  positive Borel measure. Let $\varepsilon>0$.  Then
  $\int_{-\infty}^{0}(\chi_\varepsilon\circ x)\,d\mu$ is convex. Moreover, if we
  consider the subsets of $(-\infty,0]$
  \begin{equation}
    \begin{gathered}\label{eq:subsets_aumann}
      E_\varepsilon^-=x^{-1}((-\infty,-b(\varepsilon)))\,,\quad
      E_\varepsilon^0=x^{-1}([-b(\varepsilon),b(\varepsilon)])\,,\\
      \text{and}\quad E_\varepsilon^+=x^{-1}((b(\varepsilon),\infty)),
    \end{gathered}
  \end{equation}
  then
  \[
    \int_{-\infty}^{0}(\chi_\varepsilon\circ x)\,d\mu=
    \left[(1-\varepsilon)\,\mu(E_\varepsilon^+),\,\varepsilon\,\mu(E_\varepsilon^-)+\mu(E_\varepsilon^0)+\mu(E_\varepsilon^+)\right]\,.
  \]
\end{lemma}
\begin{proof}
  Clearly, it is enough to check the second part of the statement. Let $s$ be a
  selector of $\chi_\varepsilon\circ x$. Define
  $s^-,s^+:(-\infty,0]\to\mathbb{R}$ for $t\in\mathbb{R}$ by
  \[
    s^-(t)=\left\{
      \begin{array}{ll}
        \varepsilon & \text{if } t\in E_\varepsilon^-, \\
        1 & \text{if } t\in E_\varepsilon^0\cup E_\varepsilon^+,
      \end{array}
    \right.\qquad \text{and}\qquad s^+(t)=\left\{
      \begin{array}{ll}
        0 & \text{if } t\in E_\varepsilon^-\cup E_\varepsilon^0, \\
        1-\varepsilon & \text{if } t\in E_\varepsilon^+.
      \end{array}
    \right.
  \]
  Then $s^-\leq s\leq s^+$ on $(-\infty,0]$, whence
  \[
    (1-\varepsilon)\mu(E_\varepsilon^+)=\int_{-\infty}^{0}\!\!s^-\,d\mu\leq\int_{-\infty}^{0}\!\!s\,d\mu
    \leq\int_{-\infty}^{0}\!\!s^+\,d\mu=\varepsilon\,\mu(E_\varepsilon^-)+\mu(E_\varepsilon^0)+\mu(E_\varepsilon^+).
  \]
  The first inclusion is proved. As for the converse inclusion, we consider
  three cases. First, suppose
  $\xi\in[(1-\varepsilon)\mu(E_\varepsilon^+),\mu(E_\varepsilon^+)]$. We can
  define $s:(-\infty,0]\to\mathbb{R}$ by
  \[
    s(t)=\left\{
      \begin{array}{ll}
        \frac{\xi}{\mu(E_\varepsilon^+)} & \text{if } t\in E_\varepsilon^+, \\
        0 & \text{otherwise}, 
      \end{array}
    \right.\qquad\text{or}\qquad s(t)=\left\{
      \begin{array}{ll}
        1 & \text{if } t\in E_\varepsilon^+, \\
        0 & \text{otherwise},
      \end{array}
    \right.
  \]
  depending on whether $\mu(E_\varepsilon^+)>0$ or $\mu(E_\varepsilon^+)=0$,
  resp. In any case, we have $\int_{-\infty}^{0}s\,d\mu=\xi$. On the other hand,
  the cases
  $\xi\in[\mu(E_\varepsilon^+),\mu(E_\varepsilon^0)+\mu(E_\varepsilon^+)]$ and
  $\xi\in[\mu(E_\varepsilon^0)+\mu(E_\varepsilon^+),\,\varepsilon\mu(E_\varepsilon^-)+\mu(E_\varepsilon^0)+\mu(E_\varepsilon^+)]$
  are analogous, so their proofs are omitted.  \qed
\end{proof}

\begin{lemma}\label{integral_limit_fixed_function}
  Let $x:(-\infty,0]\to\mathbb{R}$ be a measurable function and let $\mu$ be a
  positive Borel measure. Let $\{\varepsilon_m\}_m\subset[0,\frac{1}{5}]$
  decrease to $\varepsilon_0$ as $m\uparrow\infty$. Then
  \[
    \lim_{m\to \infty}\haus\left(\int_{-\infty}^{0}(\chi_{\varepsilon_m}\circ
      x)\,d\mu,\, \int_{-\infty}^{0}(\chi_{\varepsilon_0}\circ
      x)\,d\mu\right)=0\,.
  \]
\end{lemma}
\begin{proof}
  Let us check that $\mu(E^*_{\varepsilon_m})\to \mu(E^*_{\varepsilon_0})$ as
  $m\to\infty$, where $*\in\{-,0,+\}$ and $E^*_{\varepsilon_m}$ and
  $E^*_{\varepsilon_0}$ are defined as in~\eqref{eq:subsets_aumann}. Since $b$
  is strictly increasing on $\left[0,\frac{1}{5}\right]$, then
  \[
    E^-_{\varepsilon_m}\uparrow E^-_{\varepsilon_0}\,,\quad
    E^0_{\varepsilon_m}\downarrow E^0_{\varepsilon_0}\,,\quad \text{and}\quad
    E^+_{\varepsilon_m}\uparrow E^+_{\varepsilon_0} \quad \text{as
    }m\uparrow\infty\,.
  \]
  Hence, thanks to~\ref{hyp:M},
  \[
    \mu(E^-_{\varepsilon_m})\uparrow \mu(E^-_{\varepsilon_0})\,,\quad
    \mu(E^0_{\varepsilon_m})\downarrow \mu(E^0_{\varepsilon_0})\,, \quad
    \text{and}\quad
    \mu(E^+_{\varepsilon_m})\uparrow\mu(E^+_{\varepsilon_0})\quad \text{as }
    m\uparrow\infty\,.
  \]
  This fact and Lemma~\ref{convex_integral} yield the expected result.  \qed
\end{proof}

\begin{lemma}\label{integral_limit_fixed_epsilon}
  Let $\varepsilon\in[0,\frac{1}{5}]$, let $\mu$ be a positive Borel measure,
  and let $\pi:\mathbb{R}^{n^2}\to\mathbb{R}$ be a canonical projection. Suppose
  that $\{x^m\}_m\subset C_\gamma$ and $\left\|x^m-x\right\|_\gamma\downarrow 0$
  as $m\uparrow\infty$ for some $x\in C_\gamma$. Then
  \[
    \lim_{m\to \infty}\haus\left(\int_{-\infty}^{0}(\chi_{\varepsilon}\circ
      \pi\circ x^m)\,d\mu,\, \int_{-\infty}^{0}(\chi_{\varepsilon}\circ\pi\circ
      x)\,d\mu\right)=0\,.
  \]
\end{lemma}
\begin{proof}
  For each $m\in\mathbb{N}$, consider $E^-_m$, $E^0_m$, and $E^+_m$ (resp.,
  $E^-$, $E^0$, and $E^+$) defined as in~\eqref{eq:subsets_aumann} from
  $\pi\circ x^m$ (resp., $\pi\circ x$). Now, if $t\in E^0_m\cap E^-$, then
  \[
    -b(\varepsilon)-e^{-\gamma\,t}\left\|x-x^m\right\|_\gamma \leq\pi\circ
    x^m(t)-|\pi\circ x(t)-\pi\circ x^m(t)|\leq\pi\circ x(t) <-b(\varepsilon)\,.
  \]
  Hence, if
  $A_m=\{t\leq 0:-e^{\gamma\,t}\,b(\varepsilon)-\left\|x-x^m\right\|_\gamma \leq
  e^{\gamma\,t}\,(\pi\circ x)(t)<-e^{\gamma\,t}\,b(\varepsilon)\}$, we have
  $E^0_m\cap E^-\subset A_m$.

  Let us check that $A_m\downarrow \emptyset$ as $m\uparrow\infty$. Note that
  $\left\|x^m-x\right\|_\gamma\downarrow 0$ as $m\uparrow\infty$. Therefore, if
  $t\in A_{m+1}$ for some $m\in\mathbb{N}$, then
  \[
    -e^{\gamma\,t}\,b(\varepsilon)-\left\|x-x^m\right\|_\gamma \leq
    -e^{\gamma\,t}\,b(\varepsilon)-\left\|x-x^{m+1}\right\|_\gamma \leq
    e^{\gamma\,t}\,(\pi\circ x)(t)<-e^{\gamma\,t}\,b(\varepsilon)\,,
  \]
  whence $t\in A_m$. Besides, suppose that $t\in A_m$ for all
  $m\in\mathbb{N}$. Then, for all $m\in\mathbb{N}$, we have
  \[
    -e^{\gamma\,t}\,b(\varepsilon)-\left\|x-x^m\right\|_\gamma \leq
    e^{\gamma\,t}\,(\pi\circ x)(t)<-e^{\gamma\,t}\,b(\varepsilon)\,.
  \]
  Hence,
  $-e^{\gamma\,t}\,b(\varepsilon)\leq e^{\gamma\,t}\,(\pi\circ
  x)(t)<-e^{\gamma\,t}\,b(\varepsilon)$, a contradiction. As a consequence,
  $\cap_{m\in\mathbb{N}}A_m=\emptyset$. Altogether, $A_m\downarrow \emptyset$ as
  $m\uparrow\infty$, as wanted.

  Finally, since $E^0_m\cap E^-\subset A_m$ for all $m\in\mathbb{N}$, we deduce
  that $\mu(E^0_m\cap E^-)\to 0$ as $m\to\infty$. A similar argument yields
  $\mu(E^0_m\cap E^+)\to 0$, $\mu(E^+_m\cap E^-)\to 0$, and
  $\mu(E^-_m\cap E^+)\to 0$ as $m\to\infty$.

  For each $m\in\mathbb{N}$, let $s_m$ be a selector of
  $\chi_\varepsilon\circ\pi\circ x^m$. Define
  $\widetilde{s}_m:(-\infty,0]\to\mathbb{R}$ for $t\in(-\infty,0]$ by
  \[
    \widetilde{s}_m(t)=\left\{
      \begin{array}{ll}
        s_m(t) & \text{if } t\in (E^-_m\cap E^-)\cup E^0\cup(E^+_m\cap E^+)\,,\\
        0 & \text{if } t\in (E^0_m\cap E^-)\cup(E^+_m\cap E^-)\,,\\
        1 & \text{if } t\in (E^0_m\cap E^+)\cup(E^-_m\cap E^+)\,.
      \end{array}
    \right.
  \]
  Then $\widetilde{s}_m$ is a selector of $\chi_\varepsilon\circ\pi\circ
  x$. Moreover,
  \[
    \begin{split}
      \int_{-\infty}^{0}s_m\,d\mu=&\int_{-\infty}^{0}\widetilde{s}_m\,d\mu+\int_{-\infty}^{0}(s_m-
      \widetilde{s}_m)\,d\mu
      =\int_{-\infty}^{0}\widetilde{s}_m\,d\mu+\\
      &+\int_{(E^0_m\cap E^-)\cup(E^+_m\cap E^-)}s_m\,d\mu +\int_{(E^0_m\cap
        E^+)\cup(E^-_m\cap E^+)}(s_m-1)\,d\mu\,.
    \end{split}
  \]
  Notice that
  \[
    \delta_m=\int_{(E^0_m\cap E^-)\cup(E^+_m\cap E^-)}s_m\,d\mu
    +\int_{(E^0_m\cap E^+)\cup(E^-_m\cap E^+)}(s_m-1)\,d\mu\to 0
  \]
  as $m\to\infty$, thanks to the preceding argument and the boundedness of
  $s_m$. Hence,
  \[
    \haus\left(\int_{-\infty}^{0}(\chi_{\varepsilon}\circ \pi\circ x^m)\,d\mu,\,
      \int_{-\infty}^{0}(\chi_{\varepsilon}\circ\pi\circ x)\,d\mu\right)\leq
    |\delta_m|\,,
  \]
  and the result follows immediately.  \qed
\end{proof}

We consider the $\varepsilon$-inflated differential inclusion
\addtocounter{equation}{1}
\begin{equation}\label{eq:inflated}\tag*{{\upshape(\theequation)$_{\epslabel}$}}
  u'(t)\in D(u_t)+X^\varepsilon(u_t)+I(t)\,,\quad t\geq t_0\,,
\end{equation}
where $t_0\in\mathbb{R}$, $\varepsilon\geq 0$ and
$X^\varepsilon:C_\gamma\to\mathcal{P}(\mathbb{R}^{n^2})$ is given for
$u=(x,z)\in C_\gamma$ by
\[
  \begin{split}
    X^\varepsilon(u)=\prod_{i=1}^m&\Bigg(\sum_{k\neq
      i}(z_{ki}(0)-c_{ki})\,\int_{-\infty}^{0}\chi_\varepsilon
    (x_k(\tau)-\Gamma_k)\,d\mu_{ki}(\tau)\Bigg)\times \\
    &\times\prod_{\substack{i,j=1\\i\neq
        j}}^m\left\{d_{ij}\,\int_{-\infty}^{0}\sigma_\varepsilon
      (x_i(\tau)-\Gamma_i)\,d\mu_{ij}(\tau)\,[x_j(0)]^+\right\}\,.
  \end{split}
\]
We will write $X=X_0$ and consider its associated differential
inclusion~{\renewcommand{\epslabel}{0}\ref{eq:inflated}}.

The following theorem is in the line of some results in the literature (see
e.g.~\cite{paper:hk1} and~\cite{paper:hk2}).
\begin{theorem}\label{convergence}
  Let us assume hypotheses~\ref{hyp:D},~\ref{hyp:M}, and~\ref{hyp:I}. Let $T>0$, $r>0$,
  and $\{\varepsilon_m\}_m\subset\left[0,\frac{1}{5}\right]$. Suppose that
  $\varepsilon_m\downarrow\varepsilon_0$ as $m\uparrow\infty$. Let
  $\{u^m_0\}_m\subset \overline{B}_\gamma(0,r)$ be a sequence converging to
  $u_0\in C_\gamma$. Then, for any sequence $\{u^{\varepsilon_m}\}_m$ of
  solutions of the $\varepsilon_m$-inflated differential
  inclusions~{\upshape{\renewcommand{\epslabel}{\varepsilon_m}\ref{eq:inflated}}}
  with initial data $u^{\varepsilon_m}_{t_0}=u^m_0$ at $t_0\in\mathbb{R}$, there
  is a subsequence $\{u^{\varepsilon_{m_k}}\}_k$ such that
  $u^{\varepsilon_{m_k}}$ converges uniformly on $[t_0,t_0+T]$ to
  $u^*\in C([t_0,t_0+T],\mathbb{R}^{n^2})$ and $u^*$ is a solution
  of~{\upshape{\renewcommand{\epslabel}{\varepsilon_0}\ref{eq:inflated}}} on
  $[t_0,t_0+T]$ with initial datum $u^*_{t_0}=u_0$ at $t_0$.
\end{theorem}
\begin{proof}
  We write $u^{\varepsilon_m}=(x^{\varepsilon_m},z^{\varepsilon_m})$. Since
  $\left\|u^m_0\right\|_\gamma\leq r$ and $\chi_m(\mathbb{R})\subset[0,1]$ for
  all $m\in\mathbb{N}$, a straightforward adaptation of
  Proposition~\ref{general_bound} implies that there exist $k_1,k_2>0$
  independent of $\varepsilon_m$ such that
  $\left\|u^{\varepsilon_m}_t\right\|_\gamma\leq k_1\,e^{k_2(t-t_0)}$ for all
  $t\in[t_0,t_0+T]$. Consequently,
  $\left\|u^{\varepsilon_m}(t)\right\|\leq k_1\,e^{k_2\,T}$ for all
  $t\in[t_0,t_0+T]$. Let $r_1=\max\{r,\,k_1\,e^{k_2\,T}\}$.

  For each $m\in\mathbb{N}$, there is a selector $\Sigma^{\varepsilon_m}(t)$ of
  $X^{\varepsilon_m}(u^{\varepsilon_m}_t)$, $t\in[t_0,t_0+T]$, such that
  \[
    \frac{du^{\varepsilon_m}}{dt}(t)=D(u^{\varepsilon_m}_t)+\Sigma^{\varepsilon_m}(t)+I(t)\,,\quad
    t\in[t_0,t_0+T]\,.
  \]
  Since $X^{\varepsilon_m}(u^{\varepsilon_m}_t)$ is bounded by
  \[
    \max\left\{(n-1)\,(r_1+\left\|c\right\|)\,\left\|(|\mu_{ij}|)_{i\neq
          j}\right\|_\infty,\, \left\|d\right\|\,\left\|(|\mu_{ij}|)_{i\neq
          j}\right\|_\infty\,r_1\right\}
  \]
  for all $t\in[t_0,t_0+T]$, so is $\Sigma^{\varepsilon_m}$. Besides, for all
  $t\in[t_0,t_0+T]$,
  \[
    |\pi_i\circ
    D(u^{\varepsilon_m}_t)|=|-A_i(x^{\varepsilon_m}_i(t))\,x^{\varepsilon_m}_i(t)|
    \leq \sup_{s\in[-r_1,r_1]}A_i(s)\,r_1\,,
  \]
  because
  $|x^{\varepsilon_m}_i(t)|\leq \left\|u^{\varepsilon_m}(t)\right\|
  \leq\left\|u^{\varepsilon_m}_t\right\|_\gamma\leq r_1$. Analogously,
  \[
    |\pi_{ij}\circ D(u^{\varepsilon_m}_t)|\leq
    \sup_{s\in[-r_1,r_1]}B_{ij}(s)\,r_1\,.
  \]
  Finally, from~\ref{hyp:I}, it follows that
  $\sup_{t\in[t_0,t_0+T]}|I_i(t)|<\infty$ for all
  $i\in\{1,\ldots,n\}$. Altogether, there exists $K>0$ such that
  $\left\|D(u^{\varepsilon_m}_t)+\Sigma^{\varepsilon_m}(t)+I(t)\right\|\leq K$
  for all $t\in[t_0,t_0+T]$. As a result, for all $t,s\in[t_0,t_0+T]$, we have
  \[
    \left\|u^{\varepsilon_m}(t)-u^{\varepsilon_m}(s)\right\|=
    \left\|\int_{s}^{t}(D(u^{\varepsilon_m}_\tau)
      +\Sigma^{\varepsilon_m}(\tau)+I(\tau))d\tau\right\| \leq K\,|t-s|\,.
  \]
  As a consequence, $\{u^{\varepsilon_m}\}_m$ is equi-Lipschitz continuous on
  $[t_0,t_0+T]$. The same argument yields
  $\left\|\frac{du^{\varepsilon_m}}{dt}(t)\right\|\leq K$ for all
  $t\in[t_0,t_0+T]$ and all $m\in\mathbb{N}$.

  Thanks to Theorem~4 on p.13 of Aubin and Cellina~\cite{book:ac}, there exists
  a subsequence of $\{u^{\varepsilon_m}\}_m$ (denoted $\{u^{\varepsilon_m}\}_m$
  again) such that
  \[
    \begin{gathered}
      \lim_{m\to \infty}\left\|u^{\varepsilon_m}-u^*\right\|_{[t_0,t_0+T]}=0
      \quad \text{and}\\ \lim_{m\to
        \infty}\frac{du^{\varepsilon_m}}{dt}=\frac{du^*}{dt} \quad\text{weakly
        in }L^1\left([t_0,t_0+T],\mathbb{R}^{n^2}\right)\,,
    \end{gathered}
  \]
  where $u^*:[t_0,t_0+T]\to\mathbb{R}^{n^2}$ is an absolutely continuous
  function. Let us define $\left.u^*\right|_{(-\infty,t_0)}=u_0$. Then
  $\left\|u^{\varepsilon_m}-u^*\right\|_\gamma$ vanishes as $m\to\infty$, and we
  can assume that $\left\|u^{\varepsilon_m}-u^*\right\|_\gamma\downarrow 0$ as
  $m\uparrow\infty$ without loss of generality. It is noteworthy that
  $L^1([t_0,t_0+T],\mathbb{R}^{n^2})$ has the Banach-Saks property
  (cf. Szlenk~\cite{paper:s}). As a result, again taking a subsequence if
  necessary, we have
  \begin{equation}
    \label{eq:banach_saks}
    \lim_{N\to
      \infty}\frac{1}{N}\,\sum_{m=1}^{N}\frac{du^{\varepsilon_m}}{dt}=\frac{du^*}{dt}
    \quad \text{strongly in }L^1\left([t_0,t_0+T],\mathbb{R}^{n^2}\right)\,.
  \end{equation}
  Now, for each $t\in[t_0,t_0+T]$, define
  $\Sigma^*(t)=\frac{du^*}{dt}(t)-D(u^*_t)-I(t)$. We have to prove that
  $\Sigma^*(t)\in X^{\varepsilon_0}(u^*_t)$. Let us write $u^*=(x^*,z^*)$.

  Fix a canonical projection $\pi:\mathbb{R}^{n^2}\to\mathbb{R}$. First, let us
  check that
  \[
    \lim_{m\to \infty}\haus(\pi\circ X^{\varepsilon_m}(x,z),\,\pi\circ
    X^{\varepsilon_0}(x,z))=0 \quad \text{for all }(x,z)\in C_\gamma\,.
  \]
  Indeed, from Lemma~\ref{properties_hausdorff}, it suffices to check that
  \[
    \lim_{m\to
      \infty}\haus\left(\int_{-\infty}^{0}\chi_{\varepsilon_m}(x_i(\tau)-\Gamma_i)\,d\mu_{ij}(\tau),\,
      \int_{-\infty}^{0}\chi_{\varepsilon_0}(x_i(\tau)-\Gamma_i)\,d\mu_{ij}(\tau)\right)=0
  \]
  for all $i,j\in\{1,\ldots,n\}$ with $i\neq j$, but this holds thanks to
  Lemma~\ref{integral_limit_fixed_function}. In particu\-lar,
  \begin{equation}
    \label{eq:convergence_limit_fixed_u_star}
    \lim_{m\to \infty}\haus(\pi\circ
    X^{\varepsilon_m}(u^*_t),\,\pi\circ X^{\varepsilon_0}(u^*_t))=0 \quad \text{for
      all }t\in[t_0,t_0+T]\,.
  \end{equation}
  We will prove that
  $\int_{t_0}^{t_0+T}\haus(\{\pi\circ \Sigma^*(t)\},\,\pi\circ
  X^{\varepsilon_0}(u^*_t))\,dt=0$. Observe that
  \[
    \begin{split}
      \int_{t_0}^{t_0+T}\haus(\{\pi&\circ \Sigma^*(t)\},\,\pi\circ
      X^{\varepsilon_0}(u^*_t))\,dt\leq\\
      \leq & \int_{t_0}^{t_0+T}\left|\pi\circ
        \Sigma^*(t)-\frac{1}{N}\,\sum_{m=1}^{N}\pi\circ\Sigma^{\varepsilon_m}(t)\right|dt\\
      &+\int_{t_0}^{t_0+T}\haus\left(\left\{
          \frac{1}{N}\,\sum_{m=1}^{N}\pi\circ\Sigma^{\varepsilon_m}(t)
        \right\},\,\pi\circ X^{\varepsilon_0}(u^*_t)\right)\,dt\,.
    \end{split}
  \]
  From~\eqref{eq:Lipschitz_D} and \eqref{eq:banach_saks}, it follows that the
  first addend on the right hand side vanishes as $N\to\infty$. Let us check
  that the second addend also vanishes as $N\to\infty$. Fix $t\in[t_0,t_0+T]$
  and $\delta>0$. From~\eqref{eq:convergence_limit_fixed_u_star}, it follows
  that there exists $m_0\in\mathbb{N}$, which depends on $t$, such that
  $\haus(\pi\circ X^{\varepsilon_m}(u^*_t),\,\pi\circ
  X^{\varepsilon_0}(u^*_t))<\delta$ for all $m\geq m_0$. Let
  $\varepsilon=\varepsilon_{m_0}$. Clearly,
  \begin{equation}
    \label{eq:I1_I2_I3}
    \haus\left(\left\{\frac{1}{N}\,\sum_{m=1}^{N}\pi\circ\Sigma^{\varepsilon_m}(t)\right\},\,\pi\circ
      X^{\varepsilon_0}(u^*_t)\right)\leq I_1+I_2+I_3\,,
  \end{equation}
  where
  \[
    \begin{split}
      I_1=&\haus\left(\left\{\frac{1}{N}\,\sum_{m=1}^{N}\pi\circ\Sigma^{\varepsilon_m}(t)\right\},
        \,\frac{1}{N}\,\sum_{m=1}^{N}\pi\circ X^{\varepsilon}(u^{\varepsilon_m}_t)\right)\,,\\
      I_2=&\haus\left(\frac{1}{N}\,\sum_{m=1}^{N}\pi\circ
        X^{\varepsilon}(u^{\varepsilon_m}_t),\,
        \frac{1}{N}\,\sum_{m=1}^{N}\pi\circ
        X^\varepsilon(u^*_t)\right)\,,\quad\text{and}\\
      I_3=&\haus\left(\frac{1}{N}\,\sum_{m=1}^{N}\pi\circ
        X^\varepsilon(u^*_t),\, \pi\circ X^{\varepsilon_0}(u^*_t)\right)\,.
    \end{split}
  \]
  Let us check that $I_1=0$. On the one hand, thanks to
  Lemma~\ref{properties_hausdorff}, we have
  \[
    I_1\leq
    \frac{1}{N}\,\sum_{m=1}^{N}\haus\left(\left\{\pi\circ\Sigma^{\varepsilon_m}(t)\right\},
      \,\pi\circ X^{\varepsilon}(u^{\varepsilon_m}_t)\right)\,.
  \]
  On the other hand,
  $\Sigma^{\varepsilon_m}(t)\in X^{\varepsilon_m}(u^{\varepsilon_m}_t)\subset
  X^{\varepsilon}(u^{\varepsilon_m}_t)$, because
  $\chi_{\varepsilon_m}(s)\subset\chi_{\varepsilon}(s)$ for all $s\in\mathbb{R}$
  and thanks to the monotonicity of the Aumann integral with respect to
  inclusion. Summarizing, we have $I_1=0$.

  We will prove that $I_2$ vanishes as $N\to\infty$. As before,
  \[
    I_2\leq \frac{1}{N}\,\sum_{m=1}^{N}\haus\left(\pi\circ
      X^{\varepsilon}(u^{\varepsilon_m}_t),\, \pi\circ
      X^\varepsilon(u^*_t)\right)\,.
  \]
  Thus, Banach-Saks theorem guarantees that it is enough to prove that
  \[\lim_{m\to \infty}\haus \left(\pi\circ
      X^{\varepsilon}(u^{\varepsilon_m}_t),\, \pi\circ
      X^\varepsilon(u^*_t)\right)=0\,.
  \]
  From Lemma~\ref{integral_limit_fixed_epsilon}, it follows that, for all
  $i,j\in\{1,\ldots,n\}$ with $i\neq j$,
  \[
    \lim_{m\to \infty}\haus\left(\int_{-\infty}^{0}\!\chi_{\varepsilon}
      (x^{\varepsilon_m}_i(t+\tau)-\Gamma_i)d\mu_{ij}(\tau),\,
      \int_{-\infty}^{0}\!\chi_{\varepsilon}(x^*_i(t+\tau)-\Gamma_i)
      d\mu_{ij}(\tau)\right)=0\,.
  \]
  This fact, together with Proposition~\ref{properties_hausdorff}, and the
  limit $\lim_{m\to \infty}u^{\varepsilon_m}(t)=u^*(t)$, yields
  $\lim_{m\to \infty}\haus \left(\pi\circ
    X^{\varepsilon}(u^{\varepsilon_m}_t),\, \pi\circ
    X^\varepsilon(u^*_t)\right)=0$, whence $\lim_{N\to \infty}I_2=0$.

  Finally, we will check that $I_3\leq\delta$ for all $N\in\mathbb{N}$. From
  Lemma~\ref{convex_integral}, it follows that
  $\int_{-\infty}^{0}\!\chi_{\varepsilon}(x^*_i(t+\tau)-\Gamma_i)d\mu_{ij}(\tau)$
  is a convex set for all $i,j\in\{1,\ldots,n\}$ with $i\neq j$, whence
  $X^\varepsilon(u^*_t)$ is convex as well and, consequently,
  \[
    \frac{1}{N}\,\sum_{m=1}^{N}\pi\circ X^\varepsilon(u^*_t)\subset
    X^\varepsilon(u^*_t)\,.
  \]
  Therefore,
  $I_3\leq \haus(\pi\circ X^{\varepsilon}(u^*_t),\,\pi\circ
  X^{\varepsilon_0}(u^*_t))<\delta$, as desired.

  Summarizing, we have proved that, for all $t\in[t_0,t_0+T]$,
  \[
    \lim_{N\to
      \infty}\haus\left(\left\{\frac{1}{N}\,\sum_{m=1}^{N}\pi\circ\Sigma^{\varepsilon_m}(t)\right\},\,\pi\circ
      X^{\varepsilon_0}(u^*_t)\right)=0\,.
  \]
  Since $u^{\varepsilon_m}_t\in \overline{B}_\gamma(0,r_1)$ for all
  $t\in[t_0,t_0+T]$ and $\chi_{\varepsilon_m}(\mathbb{R})\subset [0,1]$ for all
  $m\in\mathbb{N}$, $\Sigma^{\varepsilon_m}$ is bounded on $[t_0,t_0+T]$ by a
  bound independent of $m$. Consequently, the convex combination
  $\frac{1}{N}\,\sum_{m=1}^{N}\pi\circ\Sigma^{\varepsilon_m}$ is bounded on
  $[t_0,t_0+T]$ by a bound independent of $N$. Similarly, since $u^*$ is bounded
  on $[t_0,t_0+T]$, from Lemma~\ref{properties_hausdorff}, it follows that
  $\pi\circ X^{\varepsilon_0}(u^*_t)$ is bounded for $t\in[t_0,t_0+T]$ by a
  bound independent of $t$.

  We are in a position to apply Lebesgue's dominated convergence
  theorem. Specifically, we obtain
  $\int_{t_0}^{t_0+T}\haus(\{\pi\circ \Sigma^*(t)\},\,\pi\circ
  X^{\varepsilon_0}(u^*_t))\,dt=0$, as wanted. Hence,
  $\haus(\{\pi\circ \Sigma^*(t)\},\,\pi\circ X^{\varepsilon_0}(u^*_t))=0$ for
  almost every $t\in[t_0,t_0+T]$. As a result,
  \[
    \Sigma^*(t)=\frac{du^*}{dt}(t)-D(u^*_t)-I(t)\in X^{\varepsilon_0}(u^*_t)
  \]
  for almost every $t\in[t_0,t_0+T]$, i.e. $u^*$ is a solution of the
  differential
  inclusion~{\upshape{\renewcommand{\epslabel}{\varepsilon_0}\ref{eq:inflated}}}. This
  completes the proof.  \qed
\end{proof}

The previous result, together with Theorem~\ref{global_existence}, allows us to
guarantee the existence of solutions of the inflated differential
inclusion~\ref{eq:inflated}.

\begin{corollary}\label{existence_heaviside}
  Assume hypotheses~\ref{hyp:D},~\ref{hyp:M}, and~\ref{hyp:I}. Let $T>0$, $r>0$,
  and $\{\varepsilon_m\}_m\subset\left[0,\frac{1}{5}\right]$. Suppose that
  $\varepsilon_m\downarrow 0$ as $m\uparrow\infty$. Let
  $\{u^m_0\}_m\subset \overline{B}_\gamma(0,r)$ be a sequence converging to
  $u_0\in C_\gamma$. Then, for any sequence $\{u^{\varepsilon_m}\}_m$ of
  solutions of the differential
  equations~{\upshape{\renewcommand{\epslabel}{\varepsilon_m}\ref{eq:sigmoidal}}}
  with initial data $u^{\varepsilon_m}_{t_0}=u^m_0$ at $t_0\in\mathbb{R}$, there
  is a subsequence $\{u^{\varepsilon_{m_k}}\}_k$ such that
  $u^{\varepsilon_{m_k}}$ converges uniformly on $[t_0,t_0+T]$ to
  $u^*\in C([t_0,t_0+T],\mathbb{R}^{n^2})$ and $u^*$ is a solution
  of~{\upshape{\renewcommand{\epslabel}{0}\ref{eq:inflated}}} on $[t_0,t_0+T]$
  with initial datum $u^*_{t_0}=u_0$ at $t_0$.
\end{corollary}
\begin{proof}
  In the light of Theorem~\ref{convergence}, it is enough to check that
  $u^{\varepsilon_m}$ is a solution of the $\varepsilon_m$-inflated
  inclusion~{\upshape{\renewcommand{\epslabel}{\varepsilon_m}\ref{eq:inflated}}}
  for all $m\in\mathbb{N}$. Fix $m\in\mathbb{N}$, $i,j\in\{1,\ldots,n\}$ with
  $i\neq j$, and $t\geq t_0$. Let us write
  $u^{\varepsilon_m}=(x^{\varepsilon_m},z^{\varepsilon_m})$. Let $E^-$, $E^0$
  and $E^+$ be defined as in~\eqref{eq:subsets_aumann} from $\varepsilon_m$ and
  the function $(-\infty,0]\to\mathbb{R}$,
  $\tau \mapsto x^{\varepsilon_m}_i(t+\tau)-\Gamma_i$. Let us define
  $f^-,f^+:(-\infty,0]\to\mathbb{R}$ for $\tau\in(-\infty,0]$ by
  \[
    f^-(\tau)=\left\{
      \begin{array}{ll}
        0 & \text{if } \tau\in E^-\cup E^0\,, \\
        1-\varepsilon_m & \text{if } \tau\in E^+\,,
      \end{array}
    \right.\quad \text{and}\quad f^+(\tau)=\left\{
      \begin{array}{ll}
        \varepsilon_m & \text{if } \tau\in E^-\,, \\
        1 & \text{if } \tau\in E^0\cup E^+\,.
      \end{array}
    \right.
  \]
  Then
  $f^-(\tau)\leq \sigma_{\varepsilon_m}(
  x^{\varepsilon_m}_i(t+\tau)-\Gamma_i)\leq f^+(\tau)$ for all $\tau\leq 0$,
  whence
  \[
    \begin{split}
      (1-\varepsilon_m)\,\mu_{ij}(E^+)=&\int_{-\infty}^{0}f^-\,d\mu_{ij}\leq
      \int_{-\infty}^{0}\sigma_{\varepsilon_m}(
      x^{\varepsilon_m}_i(t+\tau)-\Gamma_i)\,d\mu_{ij}(\tau)  \\
      \leq&\int_{-\infty}^{0}f^+\,d\mu_{ij}=\varepsilon_m\,\mu_{ij}(E^-)+\mu_{ij}(E^0)+\mu_{ij}(E^+)\,.
    \end{split}
  \]
  Consequently, from Lemma~\ref{convex_integral}, it follows that
  \[
    \int_{-\infty}^{0}\sigma_{\varepsilon_m}(
    x^{\varepsilon_m}_i(t+\tau)-\Gamma_i)\,d\mu_{ij}(\tau)\in
    \int_{-\infty}^{0}\chi_{\varepsilon_m}(
    x^{\varepsilon_m}_i(t+\tau)-\Gamma_i)\,d\mu_{ij}(\tau)\,.
  \]
  Therefore,
  $\Sigma^{\varepsilon_m}(u^{\varepsilon_m}_t)\in
  X^{\varepsilon_m}(u^{\varepsilon_m}_t)$ and the result follows immediately.
  \qed
\end{proof}

\begin{remark}
  In order to show the existence of solutions
  of~{\upshape{\renewcommand{\epslabel}{0}\ref{eq:inflated}}}, a possible
  sufficient choice in the statement of Corollary~\ref{existence_heaviside} is
  $u^m_0=u_0$ for all $m\in\mathbb{N}$.
\end{remark}

\section{Asymptotic behavior}
\label{sec:asymptotic}

We will assume hypotheses~\ref{hyp:D},~\ref{hyp:M}, and~\ref{hyp:I} throughout
this section.

\subsection{The differential inclusions}

Fix $\varepsilon\in\left[0,\frac{1}{5}\right]$. Let us consider the mapping
given by the so-called \emph{attainability set} of the $\varepsilon$-inflated
differential inclusion~\ref{eq:inflated}:
\[
  \begin{array}[t]{rccl}
    \Phi^\varepsilon: & \mathbb{R}^2_\geq\times C_\gamma & \longrightarrow & \mathcal{P}(C_\gamma) \\
                      & (t,t_0,u_0) & \mapsto         & \{v\in C_\gamma:\text{there is a solution
                                                        }u({\cdot},t_0,u_0)\text{
                                                        of~\ref{eq:inflated}}\text{ such
                                                        that}\\
                      &&&\phantom{\{v\in C_\gamma:\,}u(t_0,t_0,u_0)=u_0 \text{ and }u(t,t_0,u_0)=v\}\,.
  \end{array}
\]

Let us check that this mapping defines a set-valued process, whose properties
will be studied in the remainder of this work.

\begin{lemma}\label{set_valued_process}
  The mapping $\Phi^\varepsilon$ is a set-valued process on $C_\gamma$.
\end{lemma}
\begin{proof}
  It is trivial that $\Phi^\varepsilon(t_0,t_0,u^0)=\{u^0\}$ for all
  $u^0\in C_\gamma$ and all $t_0\in\mathbb{R}$.

  On the other hand, fix $u^0\in X$ and $t_0,t_1,t_2\in\mathbb{R}$ with
  $t_2\geq t_1\geq t_0$. Let us check that
  $\Phi^\varepsilon(t_2,t_1,\Phi(t_1,t_0,u^0))\subset
  \Phi^\varepsilon(t_2,t_0,u^0)$. Fix
  $v\in\Phi^\varepsilon(t_2,t_1,\Phi(t_1,t_0,u^0))$. There exists
  $u^2:(-\infty,t_2]\to\mathbb{R}^{n^2}$ absolutely continuous such that, for
  some selector $\Sigma_2(t)$ of $X^\varepsilon(u^2_t)$ for all $t\in[t_1,t_2]$,
  we have
  \[
    \frac{du^2}{dt}(t)=D(u^2_t)+\Sigma_2(t)+I(t)\,,\quad t\in[t_1,t_2]\,,
  \]
  $u^2_{t_1}=w\in\Phi^\varepsilon(t_1,t_0,u^0)$, and $u^2_{t_2}=v$. Similarly,
  there is a $u^1:(-\infty,t_2]\to\mathbb{R}^{n^2}$ absolutely continuous such
  that, for some selector $\Sigma_1(t)$ of $X^\varepsilon(u^1_t)$ for all
  $t\in[t_0,t_1]$, we have
  \[
    \frac{du^1}{dt}(t)=D(u^1_t)+\Sigma_1(t)+I(t)\,,\quad t\in[t_0,t_1]\,,
  \]
  $u^1_{t_0}=u^0$, and $u^1_{t_1}=w$. By concatenating $u^1$, $u^2$ and
  $\Sigma_1$, $\Sigma_2$, we conclude that
  $v\in\Phi^\varepsilon(t_2,t_0,u^0)$. The proof of the converse inclusion is
  similar.

  It only remains to check that
  $\Phi^\varepsilon(t_1,t_0,u^0)\in\mathcal{P}_c(C_\gamma)$ for all
  $(t_1,t_0)\in\mathbb{R}^2_\geq$ and all $u^0\in C_\gamma$. First, by
  Theorem~\ref{global_existence}, all the solutions of~\ref{eq:sigmoidal} are
  defined on $\mathbb{R}$. Consequently, as seen in the proof of
  Corollary~\ref{existence_heaviside}, at least some of the solutions
  of~\ref{eq:inflated} are also defined on $\mathbb{R}$. As a result,
  $\Phi^\varepsilon(t_1,t_0,u^0)\neq\emptyset$. In order to check that
  $\Phi^\varepsilon(t_1,t_0,u^0)$ is closed, fix
  $\{v^m\}_m\subset\Phi^\varepsilon(t_1,t_0,u^0)$ converging to $v\in
  C_\gamma$. Then, for each $m\in\mathbb{N}$, there is a
  $u^m:(-\infty,t_1]\to\mathbb{R}^{n^2}$ absolutely continuous such that, for
  some selector $\Sigma_m(t)$ of $X^\varepsilon(u^m_t)$ for all $t\in[t_0,t_1]$,
  we have
  \[
    \frac{du^m}{dt}(t)=D(u^m_t)+\Sigma_m(t)+I(t)\,,\quad t\in[t_0,t_1]\,,
  \]
  $u^m_{t_0}=u^0$, and $u^m_{t_1}=v^m$. By applying Theorem~\ref{convergence},
  we obtain a subsequence $\{v^{m_k}\}_k$ and a solution $u^*$
  of~\ref{eq:inflated} with initial datum $u^*_{t_0}=u^0$ and
  \[
    u^*_{t_1}=\lim_{k\to \infty}u^{m_k}_{t_1}=\lim_{k\to \infty}v^{m_k}=v
  \]
  in $C_\gamma$. Thus, $v\in\Phi^\varepsilon(t_1,t_0,u^0)$, whence
  $\Phi^\varepsilon(t_1,t_0,u^0)$ is closed, as claimed.  \qed
\end{proof}

The following lemmas are devoted to the study of the conditions required by
Theorem~\ref{abstract_attractor}. Specifically,
Lemma~\ref{upper_semi_continuous} guarantees the upper semi continuous character
of the process $\Phi^\varepsilon$, Lemma~\ref{absorbing} proves the existence of
an absorbing set for $\Phi^\varepsilon$, and
Lemma~\ref{asymptotically_upper_semi_compact} checks that $\Phi^\varepsilon$ is
asymptotically upper semi compact.

\begin{lemma}\label{upper_semi_continuous}
  The attainability set $\Phi^\varepsilon(t_1,t_0,u^0)$ is a compact subset of
  $C_\gamma$ for all $(t_1,t_0)\in\mathbb{R}^2_\geq$ and all $u^0\in
  C_\gamma$. Besides, $\Phi^\varepsilon$ is upper semi continuous.
\end{lemma}
\begin{proof}
  From Lemma~\ref{set_valued_process}, it follows that
  $\Phi^\varepsilon(t_1,t_0,u^0)$ is closed, so it is enough to check that it is
  relatively compact. Fix $\{v^m\}_m\subset\Phi^\varepsilon(t_1,t_0,u^0)$. As
  seen in the proof of Lemma~\ref{set_valued_process}, there exist a subsequence
  $\{v^{m_k}\}_k$ and a solution $u^*$ of~\ref{eq:inflated} with initial datum
  $u^*_{t_0}=u^0$ and
  $\lim_{k\to \infty}v^{m_k}=u^*_{t_1}\in\Phi^\varepsilon(t_1,t_0,u^0)$. As a
  result, $\Phi^\varepsilon(t_1,t_0,u^0)$ is compact, as expected.

  In order to check that $\Phi^\varepsilon$ is upper semi continuous, suppose
  that there \mbox{exist} $\eta>0$ and $\{u^m_0\}_m\subset C_\gamma$ converging
  to $u_0\in C_\gamma$ such that, for all $m\in\mathbb{N}$,
  $\dist_{C_\gamma}
  (\Phi^\varepsilon(t_1,t_0,u^m_0),\,\Phi^\varepsilon(t_1,t_0,u_0))\geq \eta$.
  Due to the compactness of $\Phi^\varepsilon(t_1,t_0,u^m_0)$ and
  $\Phi^\varepsilon(t_1,t_0,u_0)$, for each $m\in\mathbb{N}$, there is a
  $v^m\in \Phi^\varepsilon(t_1,t_0,u^m_0)$ such that, for all $m\in\mathbb{N}$,
  \begin{equation}
    \label{eq:upper_semi_compact_contradiction}
    \dist_{C_\gamma}(\{v^m\},\,
    \Phi^\varepsilon(t_1,t_0,u_0))=
    \dist_{C_\gamma}(\Phi^\varepsilon(t_1,t_0,u^m_0),\,
    \Phi^\varepsilon(t_1,t_0,u_0))\geq \eta.
  \end{equation}
  Notice that we can assume without loss of generality that
  $\left\|u^m_0\right\|_\gamma\leq \left\|u_0\right\|_\gamma+1$ for all
  $m\in\mathbb{N}$. As a consequence, Theorem~\ref{convergence} is still
  applicable and, again, an argument similar to that in the proof of
  Lemma~\ref{set_valued_process} implies that there are a subsequence
  $\{v^{m_k}\}_k$ and a solution $u^*$ of~\ref{eq:inflated} with initial datum
  $u^*_{t_0}=u_0$ and
  $\lim_{k\to \infty}v^{m_k}=u^*_{t_1}\in\Phi^\varepsilon(t_1,t_0,u^0)$, which
  contradicts~\eqref{eq:upper_semi_compact_contradiction}. This completes the
  proof.  \qed
\end{proof}

Let $\mathcal{M}(C_\gamma)$ be the family of nonautonomous subsets of $C_\gamma$
whose images are nonempty and closed. Henceforth, we will consider the universe
$\mathcal{D}$ of non\-autonomous subsets of $C_\gamma$ with a uniform bound on
$\mathbb{R}$, that is,
\[
  \begin{split}
    \mathcal{D}=\big\{\{D(t)\}_{t\in\mathbb{R}}\in\mathcal{M}(C_\gamma):&\,
    \text{there exists } r>0 \text{ such that }D(t)\subset
    \overline{B}_\gamma(0,r)\\
    &\,\text{for all } t\in\mathbb{R}\big\}\,.
  \end{split}
\]
Let us assume the following hypothesis, which will allow us to prove the
existence of a nonautonomous absorbing set. A weaker and simpler version of this
hypothesis will be given below.
\begin{enumerate}[label={\upshape (A)}]
\item\label{hyp:A} For each $i,j\in\{1,\ldots,n\}$ with $i\neq j$, there exist
  \begin{enumerate}[label={\upshape(\roman*)}]
  \item
    $\lambda^x_{ij},\,\lambda^z_{ij}\in\{1,\,|\mu_{ij}|^2,\,d_{ij}^2,\,|\mu_{ij}|^2\,d_{ij}^2\}$
    with $\lambda^x_{ij}\,\lambda^z_{ij}=|\mu_{ij}|^2\,d_{ij}^2$;
  \item $\nu^x_{ij},\,\nu^z_{ij}\in\{1,\,|\mu_{ij}|^2\}$ with
    $\nu^x_{ij}\,\nu^z_{ij}=|\mu_{ij}|^2$;
  \item $\eta_{ij}>0$;
  \item $\varepsilon_{ij}>0$;
  \end{enumerate}
  such that
  $ \sum_{k\neq i}(\nu^x_{ki}\,\eta_{ki}
  +\lambda^x_{ki}\,\varepsilon_{ki})<2\,\alpha$ and
  $\frac{\nu^z_{ij}}{\eta_{ij}}+\frac{\lambda^z_{ij}}{\varepsilon_{ij}}<2\,\alpha$.
\end{enumerate}

\begin{lemma}\label{absorbing}
  There exists $\gamma_0>0$ such that, for all $\gamma\in(0,\gamma_0)$, the
  set-valued process $\Phi^\varepsilon$ on $C_\gamma$ has a nonautonomous
  $\mathcal{D}$-pullback absorbing set.
\end{lemma}
\begin{proof}
  Let $u=(x,z)$ be a solution of~\ref{eq:inflated} with initial datum $u^0$ at
  $t_0\in\mathbb{R}$. Fix $i,j\in\{1,\ldots,n\}$ with $i\neq j$ and $t\geq t_0$.
  We have
  \[
    \begin{split}
      x'_i(t)\,x_i(t)\in&-A_i(x_i(t))\,x_i^2(t)
      +\pi_i\circ X^\varepsilon(u_t)\,x_i(t)+I_i(t)\,x_i(t)\,,\\
      z'_{ij}(t)\,z_{ij}(t)\in&-B_{ij}(z_{ij}(t))\,z_{ij}^2(t) +\pi_{ij}\circ
      X^\varepsilon(u_t)\,z_{ij}(t)\,.
    \end{split}
  \]
  Therefore, applying~\ref{hyp:D} and the definition of $X^\varepsilon$, we
  obtain
  \begin{equation}
    \label{eq:absorbing_bound_derivatives}
    \begin{split}
      \tfrac{1}{2}(x_i^2(t))'\leq&-\alpha\,x_i^2(t) +\sum_{k\neq
        i}(|z_{ki}(t)|+c_{ki})\,|\mu_{ki}|\,|x_i(t)|
      +|I_i(t)|\,|x_i(t)|\,,\\
      \tfrac{1}{2}(z_{ij}^2(t))'\leq&-\alpha\,z_{ij}^2(t)
      +d_{ij}\,|\mu_{ij}|\,|x_j(t)|\,|z_{ij}(t)|\,.
    \end{split}
  \end{equation}
  Thanks to~\ref{hyp:A}, there exist $\omega,\xi>0$ such that
  \begin{equation}
    \label{eq:absorbing_improved_A}
    \sum_{k\neq i}(\nu^x_{ki}\,\eta_{ki}
    +\lambda^x_{ki}\,\varepsilon_{ki})+(n-1)\,\omega+\xi<2\,\alpha\,.
  \end{equation}
  An application of Young's inequality and hypothesis~\ref{hyp:A}
  to~\eqref{eq:absorbing_bound_derivatives} yields
  \[
    \begin{split}
      (x_i^2(t))'\leq&-2\,\alpha\,x_i^2(t) +\sum_{k\neq
        i}\left[\frac{\nu^z_{ki}}{\eta_{ki}}\,z_{ki}^2(t)
        +\nu^x_{ki}\,\eta_{ki}\,x_i^2(t)+
        \frac{c_{ki}^2\,|\mu_{ki}|^2}{\omega}+\omega\,x_i^2(t)\right]\\
      &\,+ \frac{I_i^2(t)}{\xi}+\xi\,x_i^2(t)\,,\\
      (z_{ij}^2(t))'\leq&-2\,\alpha\,z_{ij}^2(t) +
      \frac{\lambda^z_{ij}}{\varepsilon_{ij}}\,z_{ij}^2(t)
      +\lambda^x_{ij}\,\varepsilon_{ij}\,x_j^2(t)\,.
    \end{split}
  \]
  By adding all these inequalities, for all $t\geq t_0$, we have
  \begin{multline*}
    \left(\left\|u(t)\right\|_2^2\right)'\leq \sum_{i=1}^{n}
    \left[-2\,\alpha+\sum_{k\neq i}(\nu^x_{ki}\,\eta_{ki}
      +\lambda^x_{ki}\,\varepsilon_{ki})+(n-1)\,\omega+\xi\right]\,x_i^2(t)\\
    +\sum_{\substack{i,j=1\\i\neq j}}^m\left[-2\,\alpha+
      \dfrac{\nu^z_{ij}}{\eta_{ij}}+\dfrac{\lambda^z_{ij}}{\varepsilon_{ij}}\right]
    \,z_{ij}^2(t) + \frac{1}{\omega}\,\max_{i=1,\ldots,n}\sum_{k\neq
      i}c_{ki}^2\,|\mu_{ki}|^2 + \frac{1}{\xi}\,\sum_{i=1}^{n}I_i^2(t)\,.
  \end{multline*}
  Hence, thanks to~\ref{hyp:A} and~\eqref{eq:absorbing_improved_A}, there is a
  $\beta>0$ such that
  \[
    \left(\left\|u(t)\right\|_2^2\right)'\leq -\beta\,\left\|u(t)\right\|_2^2
    +K+ \frac{1}{\xi}\,\sum_{i=1}^{n}I_i^2(t)\,,
  \]
  where $K= \frac{1}{\omega}\,\sum_{i=1}^{n}c_{ki}^2\,|\mu_{ki}|^2$. Now,
  integrating between $t_0$ and $t$, we obtain
  \[
    \left\|u(t)\right\|_2^2\leq e^{-\beta(t-t_0)}\,\left\|u(t_0)\right\|_2^2 +
    \frac{K}{\beta}\,\left(1-e^{-\beta(t-t_0)}\right) +
    \frac{1}{\xi}\,\sum_{i=1}^{n}\int_{t_0}^{t}e^{-\beta(t-s)}I_i^2(s)\,ds\,,
  \]
  whence, if $\gamma\in\left(0,\frac{\beta}{2}\right)$, then
  \begin{equation}\label{eq:absorbing_bound_solution}
    \left\|u(t)\right\|\leq n\,e^{-\gamma(t-t_0)}\left\|u(t_0)\right\| +
    \sqrt{\frac{K}{\beta}}+
    \frac{1}{\sqrt{\xi}}\left(\sum_{i=1}^{n}\int_{t_0}^{t}
      I_i^2(s)\,ds\right)^{\frac{1}{2}},\, t\geq t_0.
  \end{equation}
  A simple computation yields
  \begin{equation}
    \label{eq:absorbing_bound_history}
    \left\|u_t\right\|_\gamma\leq n\,e^{-\gamma(t-t_0)}\,\left\|u_{t_0}\right\|_\gamma
    + \sqrt{\frac{K}{\beta}}+
    \frac{1}{\sqrt{\xi}}\,\left(\sum_{i=1}^{n}\int_{t_0}^{t}
      I_i^2(s)\,ds\right)^{\frac{1}{2}}\,,\,t\geq t_0.
  \end{equation}
  As a consequence, the nonautonomous set defined for $t\in\mathbb{R}$ by
  \begin{equation}
    \label{eq:absorbing_set}
    \Lambda(t)=\left\{\varphi\in C_\gamma:\left\|\varphi\right\|_\gamma\leq
      \sqrt{\frac{K}{\beta}}+
      \frac{1}{\sqrt{\xi}}\,\left(\sum_{i=1}^{n}\int_{t_0}^{t}
        I_i^2(s)\,ds\right)^{\frac{1}{2}}+1\right\}
  \end{equation}
  satisfies the desired property.  \qed
\end{proof}

\begin{remark}
  One of the many particular instances of~\ref{hyp:A} is the following: for each
  $i,j\in\{1,\ldots,n\}$ with $i\neq j$, we define $\lambda^x_{ij}=1$,
  $\lambda^z_{ij}=|\mu_{ij}|^2\,d_{ij}^2$, $\nu^x_{ij}=1$,
  $\nu^z_{ij}=|\mu_{ij}|^2$, $\eta_{ij}=\varepsilon_{ij}=\eta>0$.  Then the
  following inequalities must be satisfied for~\ref{hyp:A} to hold:
  \[
    2\,(n-1)\,\eta<2\,\alpha \quad \text{and}\quad
    \frac{|\mu_{ij}|^2}{\eta}+\frac{|\mu_{ij}|^2\,d_{ij}^2}{\eta}<2\,\alpha\,,\quad
    i,j\in\{1,\ldots,n\}\,,\, i\neq j\,.
  \]
  Hence, a sufficient condition for~\ref{hyp:A} to hold is
  \[
    (n-1)\,|\mu_{ij}|^2\,(1+d_{ij}^2)<2\,\alpha^2\,,\qquad
    i,j\in\{1,\ldots,n\}\,,\quad i\neq j\,.
  \]
\end{remark}

The following lemma contains an adaptation of Ascoli-Arzel\`a's theorem to
$C_\gamma$.

\begin{lemma}\label{AA}
  Let $\{u^m\}_m$ be a sequence in $C_\gamma$ satisfying the following
  conditions:
  \begin{enumerate}[label={\upshape(\roman*)}]
  \item\label{AA_bound} there is an $r>0$ such that
    $\left\|u^m\right\|_\gamma\leq r$ for all $m\in\mathbb{N}$;
  \item\label{AA_equicontinuity} the sequence $\{e^{\gamma\,{\cdot}}\,u^m\}_m$
    is equicontinuous on $(-\infty,0]$;
  \item\label{AA_limit} the limit
    $\lim_{t\to -\infty}e^{\gamma\,t}\,u^m(t)=L_m\in\mathbb{R}^{n^2}$ is uniform
    in $m\in\mathbb{N}$, that is, for all $\eta>0$, there exists $t_0<0$ such
    that $\left\|e^{\gamma\,t}\,u^m(t)-L_m\right\|<\eta$ for all $t<t_0$ and all
    $m\in\mathbb{N}$.
  \end{enumerate}
  Then there exists a subsequence $\{u^{m_k}\}_k$ converging to $u\in C_\gamma$.
\end{lemma}
\begin{proof}
  Let $\Psi:C_\gamma\to C([-1,0],\mathbb{R}^{n^2})$ be defined for
  $u\in C_\gamma$ and $s\in[-1,0]$ by
  \[
    \Psi(u)(s)=\left\{
      \begin{array}{ll}
        e^{\gamma\,\frac{s}{1+s}}\,u\left(\frac{s}{1+s}\right) & \text{if } s\in(-1,0]\,, \\
        \lim_{t\to -\infty}e^{\gamma\,t}\,u^m(t) & \text{if } s=-1\,.
      \end{array}
    \right.
  \]
  Then $\Psi$ is an isometric isomorphism when the supremum norm,
  $\left\|{\cdot}\right\|_{[-1,0]}$, is consi\-dered on
  $C([-1,0],\mathbb{R}^{n^2})$ (cf.~\cite{book:hmn}). As a result,
  condition~\ref{AA_bound} implies that, for all $m\in\mathbb{N}$,
  $\left\|\Psi(u^m)\right\|_{[-1,0]}\leq \left\|u^m\right\|_\gamma\leq r$.

  Let us check that $\{\Psi(u^m)\}_m$ is equicontinuous on $(-1,0]$. Fix
  $s_1\in(-1,0]$ and $\eta>0$. From~\ref{AA_equicontinuity}, it follows that
  $\{e^{\gamma\,{\cdot}}\,u^m\}_m$ is equicontinuous at
  $t_1= \frac{s_1}{1+s_1}$; let $\delta>0$ be its modulus of equicontinuity for
  $\eta$ at $t_1$. Then there is a $\lambda>0$ such that, if $s_2\in(-1,0]$ and
  $|s_1-s_2|<\lambda$, then $|t_1-t_2|<\delta$, where $t_2=
  \frac{s_2}{1+s_2}$. As a consequence, for all $m\in\mathbb{N}$,
  \[
    \left\|\Psi(u^m)(s_1)-\Psi(u^m)(s_2)\right\|=
    \left\|e^{\gamma\,t_1}\,u^m(t_1)-e^{\gamma\,t_2}\,u^m(t_2)\right\|<\eta\,.
  \]
  It remains to check that $\{\Psi(u^m)\}_m$ is equicontinuous at $s_1=-1$. Fix
  $\eta>0$. By~\ref{AA_limit}, there is a $t_0<0$ such that
  $\left\|e^{\gamma\,t}\,u^m(t)-\Psi(u^m)(-1)\right\|<\eta$ for all $t<t_0$ and
  all $m\in\mathbb{N}$. Now, there is a $\lambda>0$ such that, if $s_2\in(-1,0]$
  and $|-1-s_2|<\lambda$, then $\frac{s_2}{1+s_2}<t_0$. Therefore, for all
  $m\in\mathbb{N}$,
  \[
    \left\|\Psi(u^m)(-1)-\Psi(u^m)(s_2)\right\|=
    \left\|e^{\gamma\,\frac{s_2}{1+s_2}}\,u^m\left(\tfrac{s_2}{1+s_2}\right)
      -\Psi(u^m)(s_2)\right\|<\eta\,.
  \]
  Consequently, $\{\Psi(u^m)\}_m$ is equicontinuous on $[-1,0]$.

  Finally, Ascoli-Arzel\`a's theorem guarantees the existence of
  $\{\Psi(u^{m_k})\}_{m_k}$ converging to
  $\varphi\in C([-1,0],\mathbb{R}^{n^2})$. Hence,
  $\lim_{k\to \infty}u^{m_k}=\Psi^{-1}(\varphi)$ in $C_\gamma$, as wanted.  \qed
\end{proof}

\begin{lemma}\label{asymptotically_upper_semi_compact}
  There exists $\gamma_0>0$ such that, if $\gamma\in(0,\gamma_0)$, the
  set-valued process $\Phi^\varepsilon$ on $C_\gamma$ is asymptotically upper
  semi compact.
\end{lemma}
\begin{proof}
  Fix $\{D(t)\}_{t\in\mathbb{R}}\in\mathcal{D}$, $\tau_m\to\infty$ as
  $m\to\infty$, $t\in\mathbb{R}$, and, for each $m\in\mathbb{N}$,
  $u^m_0\in D(t-\tau_m)$, $v^m\in\Phi^\varepsilon(t,t-\tau_m,u^m_0)$. Let us
  prove that $\{v^m\}_m$ has a convergent subsequence in $C_\gamma$. There
  exists $r>0$ such that $D(s)\subset \overline{B}_\gamma(0,r)$ for all
  $s\in\mathbb{R}$. It suffices to check conditions~\ref{AA_bound},
  \ref{AA_equicontinuity}, and~\ref{AA_limit} of Lemma~\ref{AA}. For each
  $m\in\mathbb{N}$, there exists $u^m:(-\infty,t]\to\mathbb{R}^{n^2}$ absolutely
  continuous such that, for some selector $\Sigma_m(s)$ of
  $X^\varepsilon(u^m_s)$ for all $s\in[t-\tau_m,t]$,
  \[
    \frac{du^m}{dt}(s)=D(u^m_s)+\Sigma_m(s)+I(s)\,,\quad s\in[t-\tau_m,t]\,,
  \]
  $u^m_{t-\tau_m}=u^m_0$, and $u^m_t=v^m$.  Let
  $r_1=n\,r+ \sqrt{\frac{K}{\beta}}+
  \frac{1}{\sqrt{\xi}}\,\left(\sum_{i=1}^{n}\int_{\mathbb{R}}
    I_i^2\right)^{\frac{1}{2}}$. As seen in~\eqref{eq:absorbing_bound_history},
  for each $m\in\mathbb{N}$,
  \[
    \left\|v^m\right\|_\gamma\leq \left\|u^m_t\right\|_\gamma\leq
    n\,e^{-\gamma\,\tau_m}\,\left\|u^m_0\right\|_\gamma+ \sqrt{\frac{K}{\beta}}+
    \frac{1}{\sqrt{\xi}}\,\left(\sum_{i=1}^{n}\int_{\mathbb{R}}
      I_i^2\right)^{\frac{1}{2}}\leq r_1\,.
  \]
  Thus, condition~\ref{AA_bound} of Lemma~\ref{AA} is met.  Now, fix $t_1\leq 0$
  and $\eta>0$. If $t_2\leq 0$, then
  \begin{equation}\label{usc_condition_ii}
    \left\|e^{\gamma\,t_1}v^m(t_1)-e^{\gamma\,t_2}v^m(t_2)\right\|=
    \left\|e^{\gamma\,t_1}u^m(t+t_1)-e^{\gamma\,t_2}u^m(t+t_2)\right\|\,.
  \end{equation}
  There are four cases. First, if $t_1,t_2\geq -\tau_m$, then,
  from~\eqref{usc_condition_ii}, it follows
  \[
    \begin{split}
      \left\|e^{\gamma\,t_1}v^m(t_1)-e^{\gamma\,t_2}v^m(t_2)\right\| \leq&\,
      e^{\gamma\,t_1}\,|1-e^{\gamma\,(t_2-t_1)}|\,\left\|u^m(t+t_1)\right\|\\
      &+e^{\gamma\,t_2}\,\left\|u^m(t+t_1)-u^m(t+t_2)\right\|\,.
    \end{split}
  \]
  
  Notice that
  $\left\|u^m(t+t_1)\right\|=\left\|v^m(t_1)\right\|\leq
  e^{-\gamma\,t_1}\,r_1$. Also, there is a $\delta>0$ such that
  $2\,e^{-\gamma\,t_1}\,r_1\,|1-e^{\gamma\,(t_2-t_1)}|\leq \eta$ if
  $|t_1-t_2|<\delta$. Besides,
  \[
    \left\|u^m(t+t_1)-u^m(t+t_2)\right\|\leq
    \left|\int_{t_1}^{t_2}\left\|D(u^m_{t+s})+\Sigma_m(t+s)+I(t+s)\right\|\,ds\right|\,.
  \]
  Fix $s$ between $t_1$ and $t_2$, and assume without loss of generality that
  $|t_1-t_2|\leq 1$. Then $s\in[t_1-1,t_1+1]\cap(-\infty,0]$ and
  \[
    \left\|u^m_{t+s}\right\|_\gamma=\left\|v^m_s\right\|_\gamma \leq
    e^{-\gamma\,s}\,\left\|v^m\right\|_\gamma\leq e^{\gamma(1-t_1)}\,r_1\,.
  \]
  A straightforward adaptation of Proposition~\ref{f_bounded} implies that there
  exists an $L>0$ independent of $\varepsilon$ such that
  \[
    \left\|D(\varphi)+X^\varepsilon(\varphi)+I(\sigma)\right\|\leq L\,, \quad
    (\sigma,\varphi)\in[t+t_1-1,t+t_1+1]\times
    \overline{B}_\gamma(0,e^{\gamma(1-t_1)}\,r_1)\,.
  \]
  As a consequence, we obtain
  $\left\|u^m(t+t_1)-u^m(t+t_2)\right\|\leq L\,|t_1-t_2|$. Therefore, if
  $|t_1-t_2|\leq\min\{\delta,\,\frac{\eta}{2\,L}\}$, then
  $\left\|e^{\gamma\,t_1}v^m(t_1)-e^{\gamma\,t_2}v^m(t_2)\right\|\leq\eta$, as
  wanted.

  Second, if $t_1,t_2\leq -\tau_m$, then, from~\eqref{usc_condition_ii}, it
  follows
  \begin{multline*}
    \left\|e^{\gamma\,t_1}v^m(t_1)-e^{\gamma\,t_2}v^m(t_2)\right\|=\\
    =e^{-\gamma\,\tau_m}\,\left\|e^{\gamma(t_1+\tau_m)}u^m_0(t_1+\tau_m)
      -e^{\gamma(t_2+\tau_m)}u^m_0(t_2+\tau_m)\right\|\leq
    2\,r\,e^{-\gamma\,\tau_m}\,.
  \end{multline*}

  Since $\lim_{m\to \infty}\tau_m=\infty$, there is an $m_0\in\mathbb{N}$ such
  that $2\,r\,e^{-\gamma\,\tau_m}<\eta$ for all $m> m_0$. On the other hand, let
  $\delta>0$ be a common modulus of uniform continuity on $(-\infty,0]$ for
  $\eta$ of $e^{\gamma\,{\cdot}}u^m_0$ for all $m\in\{1,\ldots,m_0\}$. In this
  situation, if $|t_1-t_2|<\delta$, then
  $\left\|e^{\gamma(t_1+\tau_m)}u^m_0(t_1+\tau_m)
    -e^{\gamma(t_2+\tau_m)}u^m_0(t_2+\tau_m)\right\|\leq\eta$ for all
  $m\in\{1,\ldots,m_0\}$ and, consequently,
  $\left\|e^{\gamma\,t_1}v^m(t_1)-e^{\gamma\,t_2}v^m(t_2)\right\|\leq\eta$ for
  all $m\in\mathbb{N}$, as expected.

  Finally, we can easily deal with the cases $t_2\leq -\tau_m\leq t_1$ and
  $t_1\leq -\tau_m\leq t_2$ with an application of the triangle inequality and
  of the two previous cases. As a result, condition~\ref{AA_equicontinuity} of
  Lemma~\ref{AA} is satisfied.

  Let us check that
  $\lim_{s\to -\infty}e^{\gamma\,s}\,v^m(s)=L_m\in\mathbb{R}^{n^2}$ is uniform
  in $m\in\mathbb{N}$. First, $v^m_{-\tau_m}=u^m_{t-\tau_m}=u^m_0$ for all
  $m\in\mathbb{N}$, whence
  \[
    L_m=\lim_{s\to -\infty}e^{\gamma\,s}\,v^m(s)= \lim_{s\to
      -\infty}e^{-\gamma\,\tau_m}\,e^{\gamma(s+\tau_m)}
    \,u^m_0(s+\tau_m)=e^{-\gamma\,\tau_m}\,M_m\,,
  \]
  where $M_m=\lim_{s\to -\infty}e^{\gamma\,s}\,u^m_0(s)$ exists thanks to the
  definition of $C_\gamma$. Note that $\left\|M_m\right\|\leq r$ for all
  $m\in\mathbb{N}$. Now, if $s\in[-\tau_m,0]$, then $t-\tau_m\leq t+s\leq t$ and
  an application of~\eqref{eq:absorbing_bound_solution} yields
  \[
    \begin{split}
      \|e^{\gamma\,s}\,v^m&(s)-L_m\|=
      \left\|e^{\gamma\,s}\,u^m(t+s)-e^{-\gamma\,\tau_m}\,M_m\right\|\\
      \leq&\,
      e^{\gamma\,s}\,\left\|u^m(t+s)\right\|+e^{-\gamma\,\tau_m}\,r\\
      \leq & \, e^{\gamma\,s}\left[
        n\,e^{-\gamma(s+\tau_m)}\,\left\|u^m(t-\tau_m)\right\|+
        \sqrt{\frac{K}{\beta}}+
        \frac{1}{\sqrt{\xi}}\,\left(\sum_{i=1}^{n}\int_{\mathbb{R}}
          I_i^2\right)^{\frac{1}{2}}
      \right]+e^{-\gamma\,\tau_m}\,r\\
      \leq &\, e^{\gamma\,s}\left[ n\,r+ \sqrt{\frac{K}{\beta}}+
        \frac{1}{\sqrt{\xi}}\,\left(\sum_{i=1}^{n}\int_{\mathbb{R}}
          I_i^2\right)^{\frac{1}{2}} \right]+e^{-\gamma\,\tau_m}\,r\,.
    \end{split}
  \]
  On the other hand, if $s\leq -\tau_m$, then
  \[
    \|e^{\gamma\,s}\,v^m(s)-L_m\|\leq
    e^{-\gamma\,\tau_m}\,e^{\gamma(s+\tau_m)}\,\left\|u^m_0(s+\tau_m)\right\|
    +e^{-\gamma\,\tau_m}\,r \leq 2\,e^{-\gamma\,\tau_m}\,r\,.
  \]
  Fix $\eta>0$. There is an $m_0\in\mathbb{N}$ such that
  $2\,e^{-\gamma\,\tau_m}\,r<\eta$ for all $m\geq m_0$. Besides, there exists
  $t_{m_0}<0$ such that
  $\left\|e^{\gamma\,s}v^m(s)-L_m\right\|\leq
  \frac{\eta}{2}+\frac{\eta}{2}=\eta$ for all $s<t_{m_0}$ and all $m\geq
  m_0$. Moreover, there exist $t_1,t_2,\ldots,t_{m_0-1}<0$ such that
  $\left\|e^{\gamma\,s}v^m(s)-L_m\right\|\leq\eta$ for all $s<t_m$ and all
  $m\in\{1,2,\ldots,m_0-1\}$. As a consequence, condition~\ref{AA_limit} of
  Lemma~\ref{AA} holds. An application of Lemma~\ref{AA} yields the desired
  result.  \qed
\end{proof}

Finally, we are in a position to prove the existence of a global attractor for
$\Phi^\varepsilon$.

\begin{theorem}\label{attractor_inflated}
  Under hypotheses~\ref{hyp:D},~\ref{hyp:M},~\ref{hyp:I}, and~\ref{hyp:A}, there
  exists $\gamma_0>0$ such that, if $\gamma\in(0,\gamma_0)$, the set-valued
  process $\Phi^\varepsilon$ on $C_\gamma$ has a unique $\mathcal{D}$-pullback
  attractor $\mathcal{A}=\{A(t)\}_{t\in\mathbb{R}}$ given by
  \[
    A(t)=\bigcap_{t_0\geq 0}\overline{\bigcup_{\tau\geq
        t_0}\Phi^\varepsilon(t,t-\tau,\Lambda(t-\tau))}\,, \qquad
    t\in\mathbb{R}\,,
  \]
  where $\{\Lambda(t)\}_{t\in\mathbb{R}}$ is defined
  by~\eqref{eq:absorbing_set}.
\end{theorem}
\begin{proof}
  It follows from Lemmas~\ref{upper_semi_continuous}, \ref{absorbing},
  and~\ref{asymptotically_upper_semi_compact} as well as
  Theorem~\ref{abstract_attractor}.  \qed
\end{proof}

\subsection{The sigmoidal systems}

Fix $\varepsilon\in\left(0,\frac{1}{5}\right]$ and consider the mapping
$S^\varepsilon:\mathbb{R}^2_{\geq}\times C_\gamma\to C_\gamma$,
$(t,t_0,u_0)\mapsto u_t$, where $u=u({\cdot},t_0,u_0)$ is the solution
of~\ref{eq:sigmoidal} with initial datum $u_0$ at $t_0$.

Let us check that $S^\varepsilon$ is a process (see e.g. Definition~1.1 of
Kloeden and Rasmussen~\cite{book:kr} for further details).

\begin{lemma}
  The mapping $S^\varepsilon$ is a process.
\end{lemma}
\begin{proof}
  It is clear that
  \[
    S^\varepsilon(t_0,t_0,u_0)=u_0 \quad \text{and}\quad
    S^\varepsilon(t_2,t_0,u_0)
    =S^\varepsilon(t_2,t_1,S^\varepsilon(t_1,t_0,u_0))\,,
  \]
  thanks to Theorem~\ref{sigmoidal_uniqueness}, for all
  $t_0,t_1,t_2\in\mathbb{R}$ with $t_0\leq t_1\leq t_2$ and all
  $u_0\in C_\gamma$.

  It only remains to check that $S^\varepsilon$ is continuous. Let
  $\{(t^m,t^m_0,v^m)\}_m$ be a sequence in $\mathbb{R}^2_{\geq}\times C_\gamma$
  converging to $(t,t_0,v)\in\mathbb{R}^2_{\geq}\times C_\gamma$.  Let us denote
  $u^m=u({\cdot},t^m_0,v^m)$ and $u=u({\cdot},t_0,v)$. Fix $\eta>0$. For each
  $s\leq 0$, we have
  \[
    \begin{split}
      S^\varepsilon(t^m,t^m_0,v^m)(s)=&\left\{
        \begin{array}{ll}
          v^m(t^m+s-t^m_0) & \text{if } t^m+s\leq t^m_0\,, \\
          u^m(t^m+s) & \text{if } t^m+s\geq t^m_0\,,
        \end{array}
      \right.\qquad \text{and}\\
      S^\varepsilon(t,t_0,v)(s)=&\left\{
        \begin{array}{ll}
          v(t+s-t_0) & \text{if } t+s\leq t_0\,, \\
          u(t+s) & \text{if } t+s\geq t_0\,.
        \end{array}
      \right.\\
    \end{split}
  \]
  Consequently, we can consider four different cases. First, suppose that
  $t^m+s\leq t^m_0$ and $t+s\leq t_0$. Then we have
  \[
    \begin{split}
      e^{\gamma\,s}\,\|v^m&(t^m+s-t^m_0)-v(t+s-t_0)\|\leq\\
      \leq&\, e^{-\gamma(t^m-t^m_0)}\,e^{\gamma(t^m+s-t^m_0)}\,
      \left\|v^m(t^m+s-t^m_0)-v(t^m+s-t^m_0)\right\|
      \\
      &+e^{\gamma\,s}\,\left\|v(t^m+s-t^m_0)-v(t+s-t_0)\right\|\\
      \leq &\,
      \left\|v^m-v\right\|_\gamma+e^{\gamma\,s}\,\left\|v(t^m+s-t^m_0)-v(t+s-t_0)\right\|\,.
    \end{split}
  \]
  Notice that $\lim_{m\to \infty}v^m=v$. Furthermore, let
  $M=\lim_{\sigma\to -\infty}e^{\gamma\,\sigma}\,v(\sigma)$ and fix $\rho<0$
  such that $\left\|e^{\gamma\,\sigma}\,v(\sigma)-M\right\|\leq \frac{\eta}{3}$
  for all $\sigma\leq\rho$. If $s\leq t_0-t+\rho-1$, then there is an
  $m_0\in\mathbb{N}$ such that, for all $m\geq m_0$, $t^m+s-t^m_0\leq\rho$,
  $t+s-t_0\leq\rho$, and
  $3\,\left\|M\right\|\,\left|e^{\gamma(t^m_0-t^m)}-e^{\gamma(t_0-t)}\right|<\eta$. Consequently,
  \[
    \begin{split}
      e^{\gamma\,s}&\,\|v(t^m+s-t^m_0)-v(t+s-t_0)\|
      \leq\\
      \leq&\, e^{-\gamma(t^m-t^m_0)}\,\left\|e^{\gamma(t^m+s-t^m_0)}
        v(t^m+s-t^m_0)-M\right\|
      +\left|e^{\gamma(t^m_0-t^m)}-e^{\gamma(t_0-t)}\right|\left\|M\right\|\\
      &+ e^{-\gamma(t-t_0)}\left\|M-e^{\gamma(t+s-t_0)}
        v(t+s-t_0)\right\| \leq\eta\,.
    \end{split}
  \]
  On the other hand, if $s\in[t_0-t+\rho-1,t_0-t]$, then there is an
  $m_0\in\mathbb{N}$ such that, for all $m\geq m_0$, $t^m+s-t^m_0$,
  $t+s-t_0\in[-\rho-2,0]$. Let $\delta>0$ be the modulus of uniform continuity
  of $v$ for $\eta$ on $[-\rho-2,0]$. There is an $m_1\geq m_0$ such that
  $|t^m-t|< \frac{\delta}{2}$ and $|t^m_0-t_0|< \frac{\delta}{2}$ for all
  $m\geq m_1$. Then $|(t^m+s-t^m_0)-(t+s-t_0)|<\delta$ and, as a result,
  $e^{\gamma\,s}\,\|v(t^m+s-t^m_0)-v(t+s-t_0)\|\leq\eta$, as desired.

  Second, suppose that $t^m+s\geq t^m_0$ and $t+s\geq t_0$. Then $t+s\in[t_0,t]$
  and
  \[
    \begin{split}
      e^{\gamma\,s}\,\|u^m(t^m+s)-u&(t+s)\| \leq\\
      \leq&\, \|u^m(t^m+s)-u^m(t+s)\|
      +\|u^m(t+s)-u(t+s)\|\\
      \leq&\, \|u^m(t^m+s)-u^m(t+s)\| +\|u^m-u\|_{[t_0,t]}\,.
    \end{split}
  \]
  Notice that, from Theorem~2.1 on p.43 and the remarks on p.48
  of~\cite{book:hmn}, it follows that there is an $m_0\in\mathbb{N}$ such that
  $\left\|u^m-u\right\|_{[t_0,t]}\leq \eta$ for all $m\geq m_0$. As for the
  first addend, if $t+s\geq t^m_0$, then applying
  Proposition~\ref{general_bound} with $r=\left\|u_0\right\|_\gamma+1$ and
  $T=t-t_0+2$, we conclude that there is an $r_1>0$ such that
  $\left\|u^m_\sigma\right\|_\gamma\leq r_1$ and, thanks to
  Proposition~\ref{f_bounded}, there is an $L>0$ such that
  $\left\|f_\varepsilon(\sigma,u^m_\sigma)\right\|\leq L$ for all
  $\sigma\in[t^m_0,t^m_0+T]$. Besides, there is an $m_1\geq m_0$ such that
  $L\,|t^m-t|< \eta$ for all $m\geq m_1$ and, therefore,
  \[
    \left\|u^m(t^m+s)-u^m(t+s)\right\|\leq
    \left\|\int_{t+s}^{t^m+s}f_\varepsilon(\sigma,u^m_\sigma)\,d\sigma\right\|
    \leq L\,|t^m-t|\leq\eta\,.
  \]
  On the other hand, if $t+s\leq t^m_0$, then
  \[
    \begin{split}
      \|u^m(t^m+s)-u^m&(t+s)\| \leq \\
      \leq &\,\left\|u^m(t^m+s)-u^m(t^m_0)\right\|
      +\left\|v^m(0)-v(0)\right\|\\
      &
      +\left\|v(0)-v(t+s-t^m_0)\right\|+\left\|v(t+s-t^m_0)-v^m(t+s-t^m_0)\right\|\,.
    \end{split}
  \]

  As above, there is an $m_0\in\mathbb{N}$ such that
  $\left\|u^m(t^m+s)-u^m(t^m_0)\right\|\leq\eta$ for all $m\geq m_0$. Besides,
  $t-t^m\leq t+s-t^m_0\leq 0$, so Lemma~\ref{gamma_implies_compact_open} and the
  continuity of $v$ at $0$ imply that there is an $m_1\geq m_0$ such that, for
  all $m\geq m_1$,
  \[
      \left\|v^m(0)-v(0)\right\| +\left\|v(0)-v(t+s-t^m_0)\right\|
      +\left\|v(t+s-t^m_0)-v^m(t+s-t^m_0)\right\|\leq\eta\,.
  \]
  The proofs of the cases $t^m+s\geq t^m_0$, $t+s\leq t_0$, and
  $t^m+s\leq t^m_0$, $t+s\geq t_0$ are similar, so they are omitted. The proof
  is finished.  \qed
\end{proof}

The proofs of Lemmas~\ref{upper_semi_continuous}, \ref{absorbing},
and~\ref{asymptotically_upper_semi_compact}, and
Theorem~\ref{abstract_attractor} can be easily adapted, and often simplified, to
obtain the following result on the existence of a pullback attractor for
$S^\varepsilon$.

\begin{theorem}\label{attractor_sigmoidal}
  Under hypotheses~\ref{hyp:D},~\ref{hyp:M},~\ref{hyp:I}, and~\ref{hyp:A}, there
  exists $\gamma_0>0$ such that, if $\gamma\in(0,\gamma_0)$, the process
  $S^\varepsilon$ on $C_\gamma$ has a unique $\mathcal{D}$-pullback attractor
  $\mathcal{A}=\{A(t)\}_{t\in\mathbb{R}}$ given by
  \[
    A(t)=\bigcap_{t_0\geq 0}\overline{\bigcup_{\tau\geq
        t_0}\Phi^\varepsilon(t,t-\tau,\Lambda(t-\tau))}\,, \qquad
    t\in\mathbb{R}\,,
  \]
  where $\{\Lambda(t)\}_{t\in\mathbb{R}}$ is defined
  by~\eqref{eq:absorbing_set}.
\end{theorem}

\section{Comparison of attractors}\label{sec:comparison}

Let us assume hypotheses~\ref{hyp:D},~\ref{hyp:M},~\ref{hyp:I},
and~\ref{hyp:A}. As seen in Section~\ref{sec:asymptotic},
Theo\-rem~\ref{attractor_inflated} guarantees the existence of a unique global
$\mathcal{D}$-pullback attractor for $\Phi^\varepsilon$,
$\mathcal{A}_{\Phi^\varepsilon}=\{A_{\Phi^\varepsilon}(t)\}_{t\in\mathbb{R}}$,
for each $\varepsilon\in\left[0,\frac{1}{5}\right]$. Similarly,
Theorem~\ref{attractor_sigmoidal} guarantees the existence of a unique global
$\mathcal{D}$-pullback attractor for $S^\varepsilon$,
$\mathcal{A}_{S^\varepsilon}=\{A_{S^\varepsilon}(t)\}_{t\in\mathbb{R}}$, for
each $\varepsilon\in\left(0,\frac{1}{5}\right]$. Henceforth, we will assume that
$\gamma>0$ is compatible with the existence of these attractors.

Notice also that $S^\varepsilon(t,t_0,u_0)\in\Phi^\varepsilon(t,t_0,u_0)$ for
all $\varepsilon>0$, all $(t,t_0)\in\mathbb{R}^2_\geq$, and all
$u_0\in C_\gamma$. As a result,
$A_{S^\varepsilon}(t)\subset A_{\Phi^\varepsilon}(t)$ for all $\varepsilon>0$
and all $t\in\mathbb{R}$. Analogously, it is clear that
$A_{\Phi^0}(t)\subset A_{\Phi^\varepsilon}(t)$ for all $\varepsilon>0$ and all
$t\in\mathbb{R}$.

The next result, in the line of the results presented in~\cite{paper:wky1},
establishes the convergence of the attractors for $\Phi^\varepsilon$ to those
for $\Phi^0$ as $\varepsilon\to 0^+$ in terms of the Hausdorff metric
$\Dist_{C_\gamma}$ given by
\[
  \Dist_{C_\gamma}(A,B)=\max\left\{\dist_{C_\gamma}(A,B),\,\dist_{C_\gamma}(B,A)\right\}\,,
\]
where $A$ and $B$ are nonempty and bounded subsets of $C_\gamma$.

\begin{theorem}\label{attractor_convergence}
  Under hypotheses~\ref{hyp:D},~\ref{hyp:M},~\ref{hyp:I}, and~\ref{hyp:A}, the
  attractors $\mathcal{A}_{\Phi^\varepsilon}$ converge pointwise to the
  attractor $\mathcal{A}_{\Phi^0}$ continuously in $C_\gamma$ in the Hausdorff
  metric, i.e.
  \[
    \lim_{\varepsilon\to
      0^+}\Dist_{C_\gamma}(A_{\Phi^\varepsilon}(t),\,A_{\Phi^0}(t))=0 \quad
    \text{for all }t\in\mathbb{R}\,.
  \]
\end{theorem}
\begin{proof}
  First, since $A_{\Phi^0}(t)\subset A_{\Phi^\varepsilon}(t)$ for all
  $\varepsilon>0$ and all $t\in\mathbb{R}$, we have
  \[
    \dist_{C_\gamma}(A_{\Phi^0}(t),\,A_{\Phi^\varepsilon}(t))=0\,.
  \]
  Let us check that
  $\dist_{C_\gamma}(A_{\Phi^\varepsilon}(t),\,A_{\Phi^0}(t))\to 0$ as
  $\varepsilon\to 0^+$ for all $t\in\mathbb{R}$. Suppose that there exist
  $t_0\in\mathbb{R}$, $\eta>0$, and $\varepsilon_m\downarrow 0$ as
  $m\uparrow \infty$ such that
  \[
    \dist_{C_\gamma}(A_{\Phi^{\varepsilon_m}}(t_0),\,A_{\Phi^0}(t_0))\geq\eta
  \]
  for all $m\in\mathbb{N}$. Given $m\in\mathbb{N}$, from the compactness of
  $A_{\Phi^{\varepsilon_m}}(t_0)$, it follows that there exists
  $a^m\in A_{\Phi^{\varepsilon_m}}(t_0)$ such that
  \begin{equation}
    \label{eq:attractor_convergence_contradiction}
    \dist_{C_\gamma}(\{a^m\},\,A_{\Phi^0}(t_0))=
    \dist_{C_\gamma}(A_{\Phi^{\varepsilon_m}}(t_0),\,A_{\Phi^0}(t_0))\geq\eta\,.
  \end{equation}
  As seen in Lemma~9.7 of~\cite{book:kr}, for each $m\in\mathbb{N}$, there is a
  solution $u^m:\mathbb{R}\to\mathbb{R}^{n^2}$
  of~{\renewcommand{\epslabel}{\varepsilon_m}\ref{eq:inflated}} such that
  $u^m_{t_0}=a^m$ and $u^m_t\in A_{\Phi^{\varepsilon_m}}(t)$ for all
  $t\in\mathbb{R}$. The compactness of $A_{\Phi^{\varepsilon_1}}(t_0)$, together
  with the fact that
  $a^m\in A_{\Phi^{\varepsilon_m}}(t_0)\subset A_{\Phi^{\varepsilon_1}}(t_0)$
  for all \mbox{$m\in\mathbb{N}$}, implies that there is a subsequence
  $\{a^{m_k}\}_k$ converging to $a\in A_{\Phi^{\varepsilon_1}}(t_0)$. In
  particular, $\{a^{m_k}\}_k$ is bounded. By applying Theorem~\ref{convergence}
  with $T=1$, we obtain a subsequence again denoted by $\{u^{m_k}\}_k$ which
  converges to a solution $u^*$
  of~{\renewcommand{\epslabel}{0}\ref{eq:inflated}} on $(-\infty,t_0+1]$ with
  initial datum $u^*_{t_0}=a$.

  On the other hand, as seen in the proof of Lemma~\ref{absorbing}, if
  $\{\Lambda(t)\}_{t\in\mathbb{R}}$ is the absorbing set defined
  by~\eqref{eq:absorbing_set}, then there is an $r>0$ such that
  $\Lambda(t)\subset \overline{B}_\gamma(0,r)$ for all $t\in\mathbb{R}$. Let us
  prove that $u^*_{t_0}\in\Phi^0(t_0,t_0-\tau,\overline{B}_\gamma(0,r))$ for all
  $\tau\geq 0$. Fix $\tau\geq 0$. Then
  $a^{m_k}_{-\tau}=u^{m_k}_{t_0-\tau}\in
  A_{\Phi^{\varepsilon_{m_k}}}(t_0-\tau)\subset \overline{B}_\gamma(0,r)$ for
  all $k\in\mathbb{N}$. A new application of Theorem~\ref{convergence}, yields a
  subsequence again denoted by $\{u^{m_k}\}_k$ which converges to a solution
  $\widehat{u}$ of~{\renewcommand{\epslabel}{0}\ref{eq:inflated}} on
  $(-\infty,t_0+1]$ with initial datum $u^*_{t_0-\tau}=a_{-\tau}$.

  Since $a^{m_k}_{-\tau}$ is bounded by $r$, so is $a_{-\tau}$. Hence,
  $\widehat{u}_{t_0}\in
  \Phi^0(t_0,t_0-\tau,\overline{B}_\gamma(0,r))$. Moreover, $u^{m_k}\to u^*$ on
  $[t_0,t_0+1]$ and $u^{m_k}_{t_0}=a^{m_k}\to a=u^*_{t_0}$ as $k\to\infty$,
  whence $\widehat{u}=u^*$ and, consequently,
  $u^*_{t_0}\in \Phi^0(t_0,t_0-\tau,\overline{B}_\gamma(0,r))$, as
  expected. Finally, this implies that
  \[
    \begin{split}
      \dist_{C_\gamma}(\{a\},\,A_{\Phi^0}(t_0))=&\dist_{C_\gamma}(\{u^*_{t_0}\},\,A_{\Phi^0}(t_0))\\
      \leq&\,
      \dist_{C_\gamma}(\Phi^0(t_0,t_0-\tau,\overline{B}_\gamma(0,r)),\,A_{\Phi^0}(t_0))
      \to 0 \quad \text{as }\tau\to\infty\,,
    \end{split}
  \]
  which
  contradicts~\eqref{eq:attractor_convergence_contradiction}. Consequently,
  $\lim_{\varepsilon\to
    0^+}\dist_{C_\gamma}(A_{\Phi^\varepsilon}(t),\,A_{\Phi^0}(t))=0$. The proof
  is finished.  \qed
\end{proof}

Finally, we can establish the convergence of the attractor of process
$S^\varepsilon$ to that of the set-valued process $\Phi^0$.

\begin{corollary}
  Under hypotheses~\ref{hyp:D},~\ref{hyp:M},~\ref{hyp:I}, and~\ref{hyp:A}, the
  attractors $\mathcal{A}_{S^\varepsilon}$ converge pointwise to the attractor
  $\mathcal{A}_{\Phi^0}$ upper semi continuously in $C_\gamma$, i.e.
  \[
    \lim_{\varepsilon\to
      0^+}\dist_{C_\gamma}(A_{S^\varepsilon}(t),\,A_{\Phi^0}(t))=0 \quad
    \text{for all }t\in\mathbb{R}\,.
  \]
\end{corollary}
\begin{proof}
  It suffices to notice that, for all $t\in\mathbb{R}$,
  Theorem~\ref{attractor_convergence} yields
  \[
    \begin{split}
      \dist_{C_\gamma}(A_{S^\varepsilon}(t),\,A_{\Phi^0}(t))\leq&
      \dist_{C_\gamma}(A_{S^\varepsilon}(t),\,A_{\Phi^\varepsilon}(t))
      +\dist_{C_\gamma}(A_{\Phi^\varepsilon}(t),\,A_{\Phi^0}(t))\\
      =&\,0+\dist_{C_\gamma}(A_{\Phi^\varepsilon}(t),\,A_{\Phi^0}(t)) \to 0
      \quad \text{as }\varepsilon\to 0^+\,.
    \end{split}
  \]
  This completes the proof.  \qed
\end{proof}


\begin{thebibliography}{xxx}

\bibitem{book:ac} {\scshape J.P. Aubin, A. Cellina}. {\itshape Differential
    Inclusions. Set‐Valued Maps and Viability Theory}. Grundlehren der
  Mathematischen Wissenschaften, 264. Springer‐Verlag, Berlin, 1984.
  
\bibitem{book:af} {\scshape J.P. Aubin, H. Frankowska}. {\itshape Set-Valued
    Analysis}. Systems and Control: Foundations and Applications,
  2. Birkh{\"a}user, Boston, 1990.
  
\bibitem{paper:ck} {\scshape T. Caraballo, P.E.  Kloeden}. {\itshape
    Non-autonomous attractor for integro-differential evolution equations}.
  Discrete Contin. Dyn. Syst. Ser. S 2 (2009), no. 1, 17--36.

\bibitem{paper:d} {\scshape F.S. De Blasi}. {\itshape On the differentiability
    of multifunctions}.  Pacific J. Math.  66 (1976), no. 1, 67--81.

\bibitem{book:dk} {\scshape P. Diamond, P.E. Kloeden}. {\itshape Metric spaces
    of fuzzy sets.  Theory and applications}.  World Scientific Publishing Co.,
  Inc., River Edge, NJ, 1994.

\bibitem{paper:hk1} {\scshape X. Han, P.E. Kloeden}. {\itshape Asymptotic
    behavior of a neural field lattice model with a Heaviside operator}. Phys. D
  389 (2019), 1--12.
  
\bibitem{paper:hk2} {\scshape X. Han, P.E. Kloeden}. {\itshape Sigmoidal
    approximations of Heaviside functions in neural lattice
    models}. J. Differential Equations 268 (2020), no. 9, 5283--5300.

\bibitem{paper:hskv} {\scshape A.V.M. Herz, B. Salzer, R. K{\"u}hn, J.L. van
    Hemmen}. {\itshape Hebbian learning reconsidered: representation of static
    and dynamic objects in associative neural nets}. Biol. Cybern. 60 (1989),
  457--467.
  
\bibitem{book:hmn} {\scshape Y. Hino, S. Murakami, T. Naito}, {\itshape
    Functional Differential Equations with Infinite Delay}, Lecture Notes in
  Math. 1473, Springer-Verlag, Berlin, 1991.

\bibitem{book:kr} {\scshape P.E. Kloeden, M. Rasmussen}. {\itshape Nonautonomous
    dynamical systems}.  Mathematical Surveys and Monographs, 176. American
  Mathematical Society, Providence, RI, 2011.

\bibitem{book:l} {\scshape D.S. Levine}. {\itshape Introduction to Neural and
    Cognitive Modelling}. Lawrence Erlbaum Associate, Inc., New Jersey, 1991.
  
\bibitem{paper:mp} {\scshape W.S. McCulloch, W. Pitts}. {\itshape A logical
    calculus of the ideas immanent in nervous activity}. Bull. Math. Biophys.  5
  (1943), 115--133.
  
\bibitem{paper:pv} {\scshape P. Pucci, G. Vitillaro}. {\itshape A representation
    theorem for Aumann integrals}.  J. Math. Anal. Appl.  102 (1984), no. 1,
  86--101.

\bibitem{paper:s} {\scshape W. Szlenk}. {\itshape Sur les suites faiblement
    convergentes dans l'espace L}.  (French) Studia Math.  25 (1965), 337--341.

\bibitem{paper:wky1} {\scshape X. Wang, P.E. Kloeden, M. Yang}, {\itshape
    Sigmoidal approximations of a delay neural lattice model with Heaviside
    functions}. Commun. Pure Appl. Anal. 19 (2020), 2385--2402.

\bibitem{paper:wky2} {\scshape X. Wang, P.E. Kloeden, M. Yang}, {\itshape
    Asymptotic behaviour of a neural field lattice model with
    delays}. Elec. Res. Arch. 28 (2020), 1037--1048.
  
\bibitem{book:w} {\scshape J. Wu}. {\itshape Introduction to neural dynamics and
    signal transmission delay}.  De Gruyter Series in Nonlinear Analysis and
  Applications, 6. Walter de Gruyter \& Co., Berlin, 2001.
  
\end{thebibliography}
\end{document}